\documentclass[20pt, twoside, a4paper]{article}

\ifx\TestIngCommand\undefined\relax                                     
\else\input ../../../proc/ppreamb-book.tex                               
\input init.tex \fi                                                     
\let\TestIngCommand\undefined                                             

\newtheorem{example}{Example}
\newtheorem{question}{Question}
\newtheorem{remark}{Remark}
\newtheorem{conj}{Conjecture}
\usepackage{amsmath,amscd,amssymb}
\usepackage{graphics}
\newtheorem{theo}{Theorem} 
\newtheorem{lem}{Lemma}
\newtheorem{cor}{Corollary}
\newtheorem{prop}{Proposition}
\newtheorem{defi}{Definition}

\newskip\ttglue\ttglue=.5em plus.25em minus.15em
\def\firstname#1{\def\FIRSTNAME{#1}\ignorespaces}                  
\def\lastname#1{\def\LASTNAME{#1}\ignorespaces}
\def\middleinitial#1{\def\MIDDLEINI{#1}\ignorespaces}
\def\department#1{\def\DEPARTMENT{#1}\ignorespaces}
\def\institute#1{\def\INSTITUTE{#1}\ignorespaces}
\def\address#1{\def\ADDRESS{#1}\ignorespaces}
\def\country#1{\def\COUNTRY{#1}\ignorespaces}
\def\otheraffiliation#1{\def\OTHERAFFILIATION{#1}\ignorespaces}
\def\email#1{\def\EMAIL{#1}\ignorespaces}
\newcount\autcount\autcount=0
\newcount\affcount\affcount=0
\newcount\numcount\newcount\nummcount
\newcount\nummmcount\newcount\nummmmcount
\def\writename#1#2{\ \kern-1ex\hbox{
  \csname AUthor\the#1\endcsname\
  \edef\TESTSTR{}\expandafter\ifx\csname auTHor\the#1\endcsname\TESTSTR
  \else\csname auTHor\the#1\endcsname.\ \fi 
  \csname authOR\the#1\endcsname$^{\csname AFF\the#1\endcsname}$
  \expandafter\ifx\csname corr\number#1\endcsname\relax
  \else\thanks{Corresponding author.}\ \fi
  }\ifnum#1<#2, \else\ \kern-1ex\fi}
\def\writeemail#1{
  \nummcount=0\relax\nummmcount=0\relax
  \loop\ifnum\nummcount<\autcount\advance\nummcount by1\relax
    {\expandafter\ifnum\csname AFF\the\nummcount\endcsname=#1\relax
    \global\advance\nummmcount by1\fi}\repeat
  \nummcount=0\relax\nummmmcount=0\relax
  \loop\ifnum\nummcount<\autcount\advance\nummcount by1\relax
    {\expandafter\ifnum\csname AFF\the\nummcount\endcsname=#1\relax
    \global\advance\nummmmcount by1\relax\def\blank{}\expandafter
    \ifx\csname EMAIL\the\nummcount\endcsname\blank(no e-mail)
    \else\csname EMAIL\the\nummcount\endcsname
    \fi 
    \ifnum\nummmmcount<\nummmcount; \fi\fi}\repeat}
\long\def\BeginAuthorList#1\EndAuthorList{#1\relax
  \author{\vbox{\hsize=390pt\noindent\numcount=0\relax
    \loop\ifnum\numcount<\autcount\advance\numcount by1\relax
      \writename{\numcount}{\autcount}
      \repeat}\\[2mm]
    \vbox{\small\numcount=0\relax
      \loop\ifnum\numcount<\affcount\advance\numcount by1\relax
        \vbox{{\count0=\numcount\relax
          \loop\expandafter\ifnum\csname AFF\the\count0\endcsname
            <\numcount\relax\advance\count0 by1\relax\repeat
          $^{\csname AFF\the\count0\endcsname}$}
        \def\BLANK{}\expandafter\ifx\csname DEPT\the\numcount\endcsname
          \BLANK
          \else\csname DEPT\the\numcount\endcsname, \fi
        \csname INST\the\numcount\endcsname,
        \csname ADDR\the\numcount\endcsname,
        \csname COUN\the\numcount\endcsname
        \edef\TEST{}\expandafter\ifx\csname OTHE\the\numcount\endcsname
          \TEST
          .\else;\break\csname OTHE\the\numcount\endcsname.\fi}
        \vbox{\writeemail{\numcount}}
        \repeat}\\}}
\expandafter\def\csname x1\endcsname{}
\expandafter\def\csname x2\endcsname{}
\expandafter\def\csname x3\endcsname{}
\expandafter\def\csname x4\endcsname{}
\expandafter\def\csname x5\endcsname{}
\expandafter\def\csname x6\endcsname{}
\expandafter\def\csname x7\endcsname{}
\expandafter\def\csname x8\endcsname{}
\expandafter\def\csname x9\endcsname{}
\def\Author#1#2{\global\advance\autcount by1\relax#2
  \expandafter\edef\csname AUthor\the\autcount\endcsname{\FIRSTNAME}
  \expandafter\edef\csname auTHor\the\autcount\endcsname{\MIDDLEINI}
  \expandafter\edef\csname authOR\the\autcount\endcsname{\LASTNAME}
  \expandafter\edef\csname EMAIL\the\autcount\endcsname{\EMAIL}
  \let\tempera\"\def\"{\string\"}\expandafter\ifx\csname x\DEPARTMENT
    \endcsname\relax
    \global\advance\affcount by1\relax\let\"\tempera
    \expandafter\edef\csname DEPT\the\affcount\endcsname{\DEPARTMENT}
    \expandafter\edef\csname INST\the\affcount\endcsname{\INSTITUTE}
    \expandafter\edef\csname ADDR\the\affcount\endcsname{\ADDRESS}
    \expandafter\edef\csname COUN\the\affcount\endcsname{\COUNTRY}
    \expandafter\edef\csname OTHE\the\affcount\endcsname{\OTHERAFFILIATION}
    \expandafter\edef\csname AFF\the\autcount\endcsname{\the\affcount}
  \else\expandafter\edef\csname AFF\the\autcount\endcsname{\DEPARTMENT}
  \fi\let\"\tempera\ignorespaces}
\def\CorrespondingAuthor#1#2{
  \expandafter\xdef\csname corr\number#1\endcsname{cor}
  \Author#1{#2}}
\def\PaperTitle#1{\title{\bf#1}}
\def\Category#1{\ignorespaces}
\def\keywords#1{{\noindent \emph{Keywords:}
  \def\BLANK{}\def\TEST{#1}\ifx\BLANK\TEST(n/a).\else#1\fi}}
\begin{document}
\large{                                                          
\PaperTitle{Fourier Quasicrystals on $\mathbb R^n$} 

\Category{(Pure) Mathematics}

\date{}

\BeginAuthorList
\Author1{
\firstname{Wayne}
    \lastname{Lawton}
    \middleinitial{M}
    \department{Department of the Theory of Functions, Institute of Mathematics and Computer Science}
    \institute{Siberian Federal University}
    \otheraffiliation{}
   \address{Krasnoyarsk}
    \country{Russian Federation}
    \email{wlawton@gmail.com}}
    \Author2{
     \firstname{August}
    \lastname{Tsikh}
    \middleinitial{K}
    \department{Department of the Theory of Functions, Institute of Mathematics and Computer Science}
    \institute{Siberian Federal University}
    \otheraffiliation{}
   \address{Krasnoyarsk}
    \country{Russian Federation}
    \email{atsikh@sfu-kras.ru}}
\EndAuthorList
\maketitle
\thispagestyle{empty}
\begin{abstract}
This paper has three aims. First, for $n \geq 1$ we construct a family of real-rooted 
trigonometric polynomial maps $P : \mathbb C^n \mapsto  \mathbb C^n$
whose divisors are Fourier Quasicrystals (FQ). For $n = 1$ these divisors include the first nontrivial FQ with positive integer coefficients constructed
by Kurasov and Sarnak \cite{kurasovsarnak}, and for $n > 1$ they overlap with Meyer's curved model sets \cite{meyer6} and two-dimensional \cite{meyer7} and multidimensional \cite{meyer8} crystalline measures. We prove that the divisors are FQ by directly computing their 
Fourier transforms using a formula derived in \cite{lawton}. Second, we extend the relationship between real-rootedness 
and amoebas, derived for $n = 1$ by 
Alon, Cohen and Vinzant \cite{alon1}, to the case $n > 1.$ The extension uses results in 
\cite{bushuevatsikh} about homology of complements of amoebas of algebraic sets of codimension $> 1.$
Third, we prove that the divisors of all uniformly generic real-rooted $P$
are FQ. The proof uses the formula relating Grothendieck residues and Newton polytopes derived by Gelfond and Khovanskii \cite{gelfondkhovanskii1} . Finally, we note that Olevskii and Ulanovskii \cite{olevskiiulanovskii2} have proved that all FQ with positive integer weights are divisors of real-rooted trigonometric polynomials for $n = 1$ but that the situation for $n > 1$ remains unsolved.
\end{abstract}
\noindent{\bf 2020 Mathematics Classification 52C23; 32A60; 32A27}
%
\footnote{\thanks{This work is supported by the Krasnoyarsk Mathematical Center and financed by the Ministry of Science and Higher Education of the Russian Federation (Agreement No. 075-02-2024-1429).}}
\tableofcontents
%
%
\section{Introduction}\label{sec1}
In 1974 Roger Penrose \cite{penrose}, \cite{grunbaumshephard} introduced his first aperiodic planar tilings with fivefold rotational symmetry.
In 1981 de Bruijn \cite{debruijn} described a "cut and project method" to construct Penrose tilings by projecting a subset of a cubic lattice in five dimensions onto a plane. 
In 1983 Paul Steinhardt and his student Dov Levine, inspired by these tilings, proposed a general concept of quasicrystals in a patent disclosure and published it in 1984 \cite{levine}.
Working independently, in 1984 Dan Schechtman and colleagues \cite{schechtman} synthesized substances exhibiting Bragg diffraction peaks like crystals but whose icosahedral 
symmetries were theoretically impossible for crystals. Schechtman was awarded the Chemistry Nobel Prize in 2011 for discovering these quasicrystals. In 2019 Steinhardt \cite{steinhardt} described the international expedition he led to Kamchatka in 2011 that found natural quasicrystals in meteorites. 
\\ \\
The cut and project method to construct models for both aperiodic tilings and quasicrystals are similar to the model sets studied a decade earlier by Meyer \cite{meyer1}. Also see \cite{meyer2,meyer3,moody}. The model sets are discrete but define Radon measures whose 
spectrums, defined by their Fourier transforms, are supported on dense subsets. 
Lagarias \cite{lagarias} proved that model sets are 
Besicovitch but not Bohr almost periodic. 
More recently Meyer
\cite{meyer4,meyer6,meyer7,meyer8} and Lev and Olevskii \cite{levolevskii1,levolevskii2,levolevskii3,
levolevskii4} constructed crystalline measures - those whose supports and spectrums are both discrete, and a subclass called Fourier quasicrystals. 
Favorov and colleagues
\cite{favorov1,favorov2,favorov3,favorov7,favorovkolbasina,favorovRR1,favorovRR2,kolbasina}  studied almost periodic sets and their relation to divisors of holomorphic almost periodic functions in tube domains.
\\ \\
In 2020  Kurasov and Sarnak \cite{kurasovsarnak} proved that the divisor of every real-rooted univariate trigonometric polynomial is a Fourier quasicrystal and constructed such polynomials from Lee-Yang polynomials. Shortly after in 2020 Olevskii and Ulanovskii \cite{olevskiiulanovskii2} proved that every Fourier quasicrystal on $\mathbb R$ is the divisor of a real-rooted univariate trigonometric polynomial.
In  2023 Alon, Cohen and Vinzant \cite{alon1} proved that every real-rooted exponential polynomial (called trigonometric polynomial in this paper) is the restriction of a Lee-Yang polynomial. These three results in combination completely characterize Fourier quasicrystals with positive integer weights on $\mathbb R.$
\\ \\
As summarized in the abstract this paper has three aims. Achieving these aims required using extensive notation and background knowledge. For this purpose we have included a long section titled Preliminaries. A key result is Fourier quasicrystals with positive integer coefficients are Bohr almost periodic multisets but Section 3.1 gives examples showing that the converse fails. For these examples the compactifications of multisets are not algebraic varieties. This fact supports our conjecture that the result of Olevskii and Ulanovskii referenced in the abstract holds for Fourier quasicrystals with positive integer weights on $\mathbb R^n$  for $n > 1.$ More specifically, the support of the Fourier quasicrystals are the set of commin zeroes of a set of $k$ trigonometric polynomials of on $\mathbb R^n$ where $k \geq n.$ All examples describe by Alon, Kummer, Kurasov and Vinzant in \cite{alon2} have weight $1,$ and one example has $k > n.$ In Example 3 in Subsection 3.2,  $k > n$ whenever $n > 1,$ and all weights $ = 1.$ However, in Theorem 3 in Section 5, $k = n$ and the weights equal the multiplicities of zeros of trigonometric polynomial maps $\mathbb C^n \mapsto \mathbb C^n.$
\\ \\
We strove to make the paper self contained. However, a detailed explanation of the material in Subsection \ref{subsec5.2} would increase the length of this paper to over 100 pages so we merely referenced the seminal paper of Gelfond and Khovanskii \cite{gelfondkhovanskii2}.
%
%
\section{Preliminaries}\label{sec2}
\subsection{Notation and Conventions}\label{subsec2.1}
Throughout this paper we use the following notation and facts. \\
$:=$ means 'is defined as', iff $:=$ if and only if, \\
$\emptyset :=$ empty set, $\circ :=$ composition of maps, $\simeq :=$ homeomorphism, $\cong :=$ isomorphism. \\
For sets $X, Y,$ $|X| :=$ number of elements in $X.$ $X \backslash Y := \{ x \in X: x \notin Y\}.$ \\
For a locally compact Hausdorff ($T_2$) topological space $X,$ \\
$C(X) :=$ locally convex topological vector space of continuous complex-valued functions on $X$ with topology of uniform convergence on compact subsets . This is a Fr\'echet space iff $X$ is countably compact. \\ 
$C_b(X) :=$ Banach space of bounded continuous functions with norm $||f||_\infty := \sup_{x \in X} |f(x)|.$ \\
$C_0(X) :=$ its Banach subspace of functions converging to $0$ at infinity (for every $\epsilon > 0$ there exists a compact $K_\epsilon \subset X$ such that $|f| \leq \epsilon$ on $X \backslash K_\epsilon$). \\
$C_c(X) :=$ its dense subspace of compactly supported functions, \\
$C_u(X) :=$ space of uniformly continuous $f : X \mapsto \mathbb C$ where $X$ is a uniform space (metric space or topological vector space). \\ 

$\mathcal B(X) := $ space of complex-valued Borel (finite) measures on $X,$ \\
$\mathcal R(X) := $ space of complex-valued Radon measures on $X.$ \\
$\mathbb Z, \mathbb Z_{\geq 0}, \mathbb Z_+, \mathbb Q,\mathbb R, \mathbb R_+, \mathbb C$ are integer, nonnegative integer, positive integer, rational, real, positive real, and complex numbers. \\
$\Re z, \, \Im z, |z| :=$ are the real part, imaginary part, modulus of $z \in \mathbb C.$ \\
Tuples are represented as column vectors.
$^T :=$ transpose. $m, n \in \mathbb Z_+.$ \\
$\mathbb R^n := \{\, [r_1,...,r_n]^T : r_1,...,r_n \in \mathbb R\, \}.$ \\ 
For $S \subset \mathbb R^n,$ $\dim_{\mathbb Z}(S) :=$ maximal number of $\mathbb Z$-independent elements in $S,$ \\
$\mathbb Z(S) := \mathbb Z$-module generated by $S,$ 
if $k = \dim_{\mathbb Z}(S) < \infty$ then $\mathbb Z(S) \cong \mathbb Z^k.$ \\
$L^p(\mathbb R^n) :=$ Lebesgue Banach space with norm $|| \cdot ||_p$ for $p \in [1,\infty].$ \\
$\mathbb C^m := \{\, [w_1,...,w_m]^T : w_1,...,w_m \in \mathbb C \, \}.$ \\
$\Re z := [\Re z_1,...,\Re z_n]^T, \, 
\Im z := [\Im z_1,...,\Im z_n]^T, \, z \in \mathbb C^n.$ \\
$\omega \cdot z := \sum_{k=1}^n \omega_k z_k, \ 
|z| := \sqrt {|z_1|^2+\cdots+|z_n|^2}, \ \omega, z  \in \mathbb C^n.$  \\
$B_n(x,r) :=$ open ball with center $x,$ radius $r > 0$ in $\mathbb R^n.$ \\
$B_{c,n}(z,r) :=$ open ball with center $z,$ radius $r > 0$ in $\mathbb C^n.$ 
Thus $B_{c,n}(z,r)$ is congruent isometric to $B_{2n}(0,r).$\\
$c_n := \pi^{n/2}/\Gamma(1+n/2)$ is  the Lebesgue measure of  $B_n(0,1).$
$\mathbb R^{m \times n} :=$ set of $m$ by $n$ real matrices defining maps $M : \mathbb C^n \mapsto \mathbb C^m.$ \\  
$\hbox{dim}_{\mathbb Z} M :=$ maximum number of $\mathbb Z$-independent rows of $M \in \mathbb R^{m \times n}.$ \\
$GL(m,\mathbb Z) := \{ A \in \mathbb Z^{m \times m} : \det A = \pm 1 \}$ is the general linear group over $\mathbb Z.$ \\
$I_m \in GL(m,\mathbb Z)$ is the identity matrix. \\
$C^* := \mathbb C \backslash \{0\}$ and $C^{*m} := (\mathbb C^*)^m$ are topological groups under multiplication. \\ 
$\mathbb T := \{z \in \mathbb C^* : |z| = 1\}$ and $\mathbb T^m$ are compact subgroups. \\
$1$ is the identity and multiplication is denoted by juxtaposition. \\
$\rho_m : \mathbb C^m \mapsto \mathbb C^{*m}$ is the homomorphism with kernel $\mathbb Z^m$ defined by
\begin{equation}\label{rhom}
	\rho_m(w) := [e^{2\pi i w_1},...,e^{2\pi i w_m}]^T, \ \ w \in \mathbb C^m.
\end{equation}
$N \in \mathbb Z^{m \times n}$ defines a homomorphism
$\widetilde N : \mathbb C^{*m} \mapsto \mathbb C^{*n}$ defined by
\begin{equation}\label{tildeN}
	(\widetilde N z)_j := \prod_{k=1}^m z_k^{N_{k,j}}, \ \ j = 1,...,n.
\end{equation}
Thus $\ell \in \mathbb Z^m$ defines the homomorphism 
$\widetilde \ell : \mathbb C^{*m} \mapsto \mathbb C^*$ 
by  
$$\widetilde \ell(z) = z^\ell := z_1^{\ell_1}\cdots z_m^{\ell_m}, \ \ z \in \mathbb C^{*m}.$$
Furthemore
\begin{equation}\label{rhomtildeN}
	\rho_m \circ N = \widetilde {N^T} \circ \rho_n, \ \ N \in Z^{m \times n}.
\end{equation}
$\mathcal L_m :=$ algebra of Laurent polynomials consisting of sums
\begin{equation}\label{q}
q = \sum_{\ell \in \Omega(q)} c_\ell \, \widetilde \ell
\end{equation}
where the spectrum $\Omega(q) \subset \mathbb Z^m$ is finite
and $c_\ell \in \mathbb C^*$ is the amplitude corresponding to frequency $\ell.$ \\
$\mathcal N(q) := $ convex hull of $\Omega(q)$ is the Newton polytope of $q.$ \\
$\mathcal V(q) := $set of vertices of $\mathcal N(q).$ \\
$\Lambda(q) := \{\lambda \in \mathbb C^{*m}
 : q(\lambda) = 0\}$ is the set of roots of $q.$ \\
$\mathcal L_m^n :=$ set of Laurent maps $Q = [q_1,...,q_n]^T  : \mathbb C^{*m} \mapsto \mathbb C^n.$ \\
$\Omega(Q) := \cup_{j=1}^n \Omega(q_j), \ 
\mathcal N(Q) := [\mathcal N(q_1),...,\mathcal N(q_n)]^T,$ \\ 
$\Lambda(Q) := \cap_{j=1}^n \Lambda(q_j)$ is the set of roots of $Q.$ \\ \\
$\omega \in \mathbb R^n$ induces a homomorphism
$\widehat \omega : \mathbb C^n \rightarrow \mathbb C^*$ defined by
\begin{equation}\label{hatw1}
	\widehat \omega (z) = e^{2\pi i\, \omega \cdot z}, \ \ z \in C^n,
\end{equation}
thus
\begin{equation}\label{hatw2}
	\widetilde \ell \circ \rho_m \circ M = \widehat {M^T\ell}, \ \ 
\ell \in \mathbb Z^m, \ M \in \mathbb R^{m \times n}.
\end{equation}
$\mathcal T_n :=$ algebra of trigonometric polynomials $p  : \mathbb C^n \mapsto \mathbb C$ consisting of sums\\
\begin{equation}\label{p}
p = \sum_{\omega \in \Omega(p)} c_\omega \, \widehat \omega
\end{equation}
where the spectrum $\Omega(p) \subset \mathbb R^m$ is finite
and $c_\omega \in \mathbb C^*$ is the amplitude corresponding to frequency $\omega.$ \\
$\mathcal N(p) := $ convex hull of $\Omega(p)$ is the Newton polytope of $p.$ \\
$\mathcal V(p) := $set of vertices of $\mathcal N(p).$ \\
$\Lambda(p) := \{\lambda \in \mathbb C^n : p(\lambda) = 0\}$ is the set of roots of $p.$ \\
$\mathcal T_n^n :=$ set of trigonometric maps $P = [p_1,...,p_n]^T  : \mathbb C^n \mapsto \mathbb C^n.$ \\
$\Omega(P) := \cup_{j=1}^n \Omega(p_j), \ 
\mathcal N(P) := [\mathcal N(p_1),...,\mathcal N(p_n)]^T,$ \\ 
$\Lambda(P) := \cap_{j=1}^n \Lambda(p_j)$ is the set of roots of $P.$ \\ \\
$P$ is real-rooted if $\Lambda(P) \subset \mathbb R^n.$ 
If $P$ is real-rooted and $\Lambda(P)$ is discrete then \\
$div\, P := \sum_{\lambda \in \Lambda} m_\lambda(P) \delta_\lambda \in \mathcal R(\mathbb R^n)$ is the divisor of $P.$ Here $\delta_\lambda$ is the point measure at $\lambda$ defined in \ref{subsec2.2} and
$m_\lambda(P)$ is the multiplicity of $\lambda$ defined in \ref{subsec2.7}.
\subsection{Tempered Distributions and Radon Measures}\label{subsec2.2}
For a topological vector space $V$ \\
$V^* :=$ its dual space consisting of continuous linear functionals $f : V \mapsto \mathbb C.$ \\
$<f,h> \, := f(h), \ \ f \in V^*, \ h \in V.$ \\
$\mathcal S(\mathbb R^n) :=$ space of Schwartz \cite{schwartz} functions $h : \mathbb R^n \mapsto \mathbb C$ all of whose derivatives exist and decay fast. Its topology is defined by seminorms
\begin{equation}
||h||_N := \sup_{x \in \mathbb R^n} 
(1+|x|)^N) \sum_{\ell_1+ \cdots + \ell_n \leq N} |D^\ell h(x)|, \ N \geq 0.
\end{equation}
$\mathcal S(\mathbb R)^*  :=$ space of tempered distributions.
Identify $S(\mathbb R^n) \subset S(\mathbb R^n)^*$ by
\begin{equation}
h \mapsto \int_{\mathbb R^n} f(x) \, h(x)\, dx, \ f, h \in S(\mathbb R^n).
\end{equation}
$\mathcal F : \mathcal S(\mathbb R^n) \mapsto \mathcal S(\mathbb R^n)$ is the Fourier transform defined by 
\begin{equation}
	\mathcal F(h)(y) := \int_{\mathbb R^n} h(x) \, e^{-2\pi i x \cdot y}\, dx, \ \ h \in \mathcal S(\mathbb R^n),\, y \in \mathbb R^n.
\end{equation}
It is a continuous bijection and $\mathcal F^{-1}(h)(x) = \int_{\mathbb R^n} h(y) \, e^{2\pi i x \cdot y} dx.$ \\
For $f, h \in \mathcal S(\mathbb R^n)$ Fubini's theorem (\cite{rudin2}, p. 164) gives
$<\mathcal F(f),h> \, :=$
$$\int_{\mathbb R^n} \int_{\mathbb R^n} e^{-2\pi i x \cdot y} f(x)h(y)dx dy = 
\int_{\mathbb R^n} \int_{\mathbb R^n} e^{-2\pi i x \cdot y} h(y)f(x)dy dc := \, <f,\mathcal F(h)>.$$
This allows us to consistently define the bijection $\mathcal F : \mathcal S(\mathbb R^n)^* \mapsto \mathcal S(\mathbb R^n)^*$ by
\begin{equation}\label{duality2}
<\mathcal F(f),h> \, := \, <f,\mathcal F(h)>, \ \ f \in \mathcal S(\mathbb R^n)^*, \, h \in \mathcal S(\mathbb R^n).
\end{equation}
Let $\widehat f := \mathcal F(f)$ for 
$f \in \mathcal S(\mathbb R^n)^*$ \\
$\mathcal S_c(\mathbb R^n) := C_c(\mathbb R^n) \cap \mathcal S(\mathbb R^n)$ is a dense subspace of $\mathcal S(\mathbb R^n).$ The 
Fourier-Laplace transform $\widehat h : \mathbb C^n \mapsto \mathbb C$ of $h \in \mathcal S_c(\mathbb R^n),$ defined by
\begin{equation}\label{FLaplace}
	\widehat h(z) := \int_{x \in \mathbb R^n} h(x)\, e^{-2\pi i \, x \cdot z}\, dx, \ \ z \in \mathbb C^n,
\end{equation}
is the unique analytic continuation of $\widehat h$ to an entire function.
\begin{lem}\label{PWS} 
$h \in \mathcal S_c(\mathbb R^n)$ implies there exists $\gamma > 0$ depending only on the support of $h,$ and for every
$N > 0$ there exists $C_N  > 0$ such that
\begin{equation}\label{pwsin}
	|\widehat h(z)| \leq C_N(1+|z|)^{-N} e^{\gamma\, ||\Im z||}, \ \ z \in \mathbb C^n.
\end{equation}
Conversely, (\ref{pwsin}) implies that $h \in \mathcal S_c(\mathbb R^n).$
\end{lem}
Proof. This is the Paley--Wiener-Schwartz theorem (\cite{hormander}, Theorem 7.3.1).
\\ \\
If $\mu \in \mathcal B(X)$ then its total variation $|\mu|,$ defined in (\cite{rudin2}, p. 116),
is a positive Borel measure, $|\mu|(X) < \infty,$ and $\mathcal B(X)$ is a Banach space with norm 
$||\mu|| := |\mu|(X)$ (\cite{rudin2}, Theorem 6.2). 
The Riesz representation theorem (\cite{rudin2}, Theorem 6.19) implies 
that $\mathcal B(X) = C_0(X)^*.$
The point measure at $x \in X$
$$\delta_x(S) := \begin{cases}
 1, \ \ x \in S \\
0, \  \ x\notin S
\end{cases}
$$
\\
for every Borel $S \subset X.$ Therefore $<\delta_x, h> = h(x), \ h \in C_0(X).$
\\
$\mathcal R(X) = C_c(X)^*$ where $C_c(X)$ has the inductive limit topology (\cite{treves}, p.515). A Radon measure is discrete if it
equals $\sum_{x \in D} c_x \, \delta_x$
where $D \subset X$ is discrete (has no limit points)
and $c_x \in \mathbb C^*$ are called its coefficients or weights. Lebesgue measure 
on $\mathbb R^n$ is a Radon measure but not a Borel measure. 
The following is from (\cite{baake}, Definition 2.1).
\begin{defi}\label{trm1}
$\mu \in \mathcal R(\mathbb R^n)$ is a tempered Radon measure if
there exists \\
$T \in \mathcal S(\mathbb R^n)^*$ such that
$<\mu, h> \, = \, <T, h>, \ \ h \in \mathcal S_c(\mathbb R^n).$
\end{defi}
Since $\mathcal S_c(\mathbb R^n)$ is dense in $\mathcal S(\mathbb R^n),$ tempered Radon measures are tempered distributions. H\"{o}rmander (\cite{hormander}, Theorem 2.1.7)  proved that compactly supported positive distributions are positive Borel measures, hence every positive tempered distribution is a positive tempered Radon measure.
\begin{lem}\label{trm2} 
For positive $\mu \in \mathcal R(\mathbb R^n)$
the following conditions are equivalent.
\begin{enumerate} 
\item $\mu$ is a tempered Radon measure.
\item There exist $C, \beta > 0$ such that
$\mu(B_n(0,r)) < C(1+r^\beta), \ \ r > 0.$
\item There exists $\alpha > 0$ such that
$(1 + |x|)^{-\alpha} \mu \in \mathcal B(\mathbb R^n).$
\end{enumerate}
\end{lem}
Favorov (\cite{favorov5}, Lemma 1) called property 2 slowly increasing and
noted its equivalence to property 1. The equivalence of properties 2 and 3 follows from 
$$
((1 + |x|)^{-\alpha} \mu)\, (\mathbb R^n) = 
\int_{r = 0}^\infty (1 + r)^{-\alpha} d \mu(B_n(0,r)).
$$
\begin{remark}\label{vmt1}
Meyer (\cite{meyer4}, p. 1) observed that  
$\mu := \sum_{k = 1}^\infty  2^k\, (\delta_{k-2^{-k}} - \delta_k)$
is a tempered Radon measure but $|\mu|$ is not because
$|\mu|(B_n(0,r)) \approx 2^{r+1}$ which violates  condition 3 in Lemma \ref{trm2}.
\end{remark}
If $\mu \in \mathcal R(\mathbb R^n)$ is Lebesgue measure then the Fourier inversion formula implies
$$<\widehat \mu, h> \, := \, <\mu, \widehat h> \, = h(0), \ \ h \in \mathcal S(\mathbb R^n)$$
hence $\widehat \mu = \delta_0.$ The Dirac comb $\mu := \sum_{\ell \in \mathbb Z^n} \delta_\ell$ is clearly a Radon measure.
\begin{lem}\label{PSF1} 
If $\mu$ is the Dirac comb then $\widehat \mu = \mu,$
or equivalently
\begin{equation}\label{PSF2}
\sum_{\ell \in \mathbb Z^n} \widehat h(\ell) = \sum_{\ell \in \mathbb Z^n} h(\ell), \ \ h \in  \mathcal S(\mathbb R^n).
\end{equation}
(\ref{PSF2}) is called the Poisson Summation Formula.
\end{lem}
Proof. Lemma \ref{PSF1} is proved in (\cite{hormander}, Theorem 7.2.1). (\ref{PSF2}) then follows since \\
$\sum_{\ell \in \mathbb Z^n} \widehat h(\ell) := \, <\mu, \widehat h> \, := \, <\widehat \mu, h> \, =  \,
<\mu,h> \, :=  \sum_{\ell \in \mathbb Z^n} h(\ell).$

\subsection{Crystalline Measures and Fourier Quasicrystals}\label{subsec2.3}
\begin{defi}\label{cmfqRn}
A Crystalline Measure (CM) on $\mathbb R^n$ is a discrete tempered Radon measure whose Fourier transform is 
a discrete Radon measure. A Fourier Quasicrystal (FQ) on $\mathbb R^n$ is a CM $\mu$
on $\mathbb R^n$ such that both $|\mu|$ and $|\widehat \mu|$ are tempered. 
\end{defi}
\begin{remark}\label{fav1}
Favorov \cite{favorov6} constructed a CM on $\mathbb R$ that is not a FQ.
\end{remark}
Lemma \ref{PSF1} implies that the Dirac comb is a FQ.
For $M \in GL(n,\mathbb R)$ and $y \in \mathbb R^n$ the affine transformation $x \mapsto Mx+y$ on $\mathbb R^n$ gives am affine transformation on discrete Radon measures
$\mu = \sum_{x\in D} c_x \delta_x \mapsto \mu_1 := \sum_{x \in D} c_x \, \delta_{Mx + y}.$
If $\mu$ is a CM with 
$\widehat \mu = \zeta := \sum_{s \in S} a_s \delta_s$ 
then a computation gives 
\begin{equation}\label{PSF3}
\widehat \mu_1 = \zeta_1 := \frac{1}{|\det M|} \sum_{s \in M^{-T}S} a_s\, e^{-2\pi i s \cdot y} \delta_s
\end{equation}
so $\mu_1$ is a CM. Furthermore if $\mu$ is a FQ then $\mu_1$ is a FQ.
Clearly finite linear combinations of CMs (resp. FQs) are CMs (resp. FQs). 
Furthermore, if $\mu$ is a CM (resp. FQ) on $\mathbb R^n$ and $\zeta$ is a CM (resp. FQ) on $\mathbb R^m,$ 
then their product $\mu \times \zeta$ is a CM (resp. FQ) on $\mathbb R^{n+m}.$
\begin{defi}\label{trivial}
A CM on $\mathbb R$ is called trivial if it is a finite linear combination of affine transformations of Dirac combs multiplied by trigonometric polynomials. Then its support is a finite union of arithmetic sequences and the CM is  a FQ. A CM on $\mathbb R^n$ for $n > 1$ is called trivial if it is a finite union of affine transformations of multidimensional Dirac combs multiplied by trigonometric polynomials; or a finite linear combination of of affine transformation of products of (not necessarily trivial) CMs on lower dimensional spaces. A FQ is called trivial if it is a trivial CM. It is called nontrivial if it is not trivial.
\end{defi}
\subsection{Bohr and Besicovitch Almost Periodicity}\label{subsec2.4}
The theory of uniform almost periodic functions on $\mathbb R$ was pioneered in three papers
\cite{bohr1} by Harald Bohr and extended to mean-squared almost periodic functions in \cite{besicovitch} by Abram Besicovitch. We summarize their results extended to 
$\mathbb R^n,$ define almost periodic measures and relate them to tempered distributions and Fourier quasicrystals in Propositions \ref{prop1} and \ref{prop2}.  
\begin{defi}\label{BAPdef}
$S \subset \mathbb R^n$ is relatively dense if there exists $r > 0$ such that 
$\mathbb R^n = \bigcup_{\tau \in S} B_n(\tau,r).$ 
$f \in C(\mathbb R^n)$ is almost periodic (\cite{bohr2}, p. 32, Main Definition) if for every $\epsilon > 0$ the set
\begin{equation}\label{translates}
	S_f(\epsilon) := \{\tau \in \mathbb R^n : ||f(\cdot-\tau) - f||_\infty \leq \epsilon\}
\end{equation}
is relatively dense. $S_f(\epsilon)$ is called the set of $\epsilon$-almost periods of $f.$
\end{defi}
We call these Bohr almost periodic functions, denote the set of them 
by $B(\mathbb R^n),$ and record that
\begin{enumerate}
\item (\cite{bohr2}, Theorem I) $B(\mathbb R^n) \subset C_b(\mathbb R^n).$
\item (\cite{bohr2},Theorem II) $B(\mathbb R^n) \subset C_u(\mathbb R^n).$
\item (\cite{bohr2},Theorem III) $B(\mathbb R^n)$ is  closed under addition, so 
$\mathcal T_n(\mathbb R^n) \subset B(\mathbb R^n),$
\item (\cite{bohr2},Theorem IV) $B(\mathbb R^n)$ is closed under multiplication
\item (\cite{bohr2},Theorem V) $B(\mathbb R^n)$ is closed under uniform 
closure.
\item (\cite{bohr2}, p. 39, 44) every $ f \in B(\mathbb R^n)$ has a mean value
\begin{equation}\label{mean}
\mathcal M(f) := \lim_{r \rightarrow \infty} \frac{r^{-n}}{c_n} \int_{x \in B_n(y,\, r)} f(x) \, dx, \ \ y \in \mathbb R^n,
\end{equation} 
where the value is indepedent of $y \in \mathbb R^n$ and convergence is uniform with respect to $y.$
\item (\cite{bohr2}, p. 77-80) The Fourier-Bohr transform $f \mapsto \mathcal F_B(f) : \mathbb R^n \mapsto \mathbb C$
\begin{equation}\label{FBf}
\mathcal F_B(f)(\omega) := \mathcal M(f(x)e^{-2\pi i\, \omega \cdot x}),  \ \ f \in B(\mathbb R^n\ \omega \in \mathbb R^n)
\end{equation}
is injective, therefore $\mathcal F_B(f)$  uniquely determines $f.$
\item (\cite{bohr2}, p. 60-65) This uniqueness is equivalent to Parseval's equation
\begin{equation}\label{Parseval}
	\mathcal M(|f|^2) = \sum_{\omega \in \Omega(f)} |\mathcal F_B(f)(\omega)|^2, \ \ f \in B(\mathbb R^n).
\end{equation}
\item (\cite{bohr2}, p. 80-88, Fundamental Theorem) $B(\mathbb R^n) = $ uniform closure of $\mathcal T_n.$ \\
\end{enumerate}
(\ref{Parseval}) implies that the spectrum 
\begin{equation}\label{spectrum}
	\Omega(f) := 
\{\omega \in \mathbb R^n : \mathcal F_B(f)(\omega) \neq 0\}
\end{equation}
is countable. 
For $f \in C_b(\mathbb R^n)$ we call $O(f) :=$ the set of translates of $f$ its orbit and
the uniform closure $\overline {O(f)}$ its orbit closure.
Bochner \cite{bochner} proved 
\begin{equation}\label{Bochner}
B(\mathbb R^n) = \{f \in C(\mathbb R^n) : \overline {O(f)} \hbox{ is compact.} \}.
\end{equation}
Besicovitch \cite{besicovitch}  defined generalized almost periodic functions
$Bes(\mathbb R^n) :=$ completion of $\mathcal T_n$ with respect to the norm
\begin{equation}
||f||_{Bes} := \sqrt {\mathcal M(|f|^2)}.
\end{equation}
Clearly $B(\mathbb R^n) \subset Bes(\mathbb R^n).$ Lagarias \cite{lagarias} proved that functions associated with
quasicrystals and Meyer's model sets are in $Bes(\mathbb R^n) \backslash B(\mathbb R^n).$ Meyer \cite{meyer3} explains this result iin detail.
\begin{remark} 
A subset $S$ of an abelian topological group $T$ is syndetic if there exists a compact $K \subset T$ such that $S+K = T.$
For $T = \mathbb R^n$ syndetic means relatively dense. These concepts play a central role in topological dynamics where thier relation to minimal systems was pioneered by Gottschalk and Hedlund (\cite{gottschalk}, Chapter 4) .
For $x \in \mathbb R^n$ and $f \in C_b(\mathbb R^n)$ define the translation of $f$ by $x,$
$\pi^xf := f(\cdot - x),$
and the orbit of $f$, 
$O(f) := \{\pi^xf : x \in \mathbb R^n\}.$
Then $\{\pi^x : x \in \mathbb R^n\}$ is a uniformly equicontinuous group of transformations of $C_b(\mathbb R^n)$ (\cite{gottschalk}, 11.08) so for every $x \in \mathbb R^n$ and every $f \in C_b(\mathbb R^n),$
$\pi^x$ extends to a map $\pi^x : \overline {O(f)} \mapsto 
\overline {O(f)}.$
(\cite{gottschalk}, 4.44) implies that $f \in B(\mathbb R^n)$ iff $\overline {X_f}$ is compact. This gives Bochner's
\cite{bochner} characterization of almost periodicity: $f$ is Bohr almost periodic if every sequence of translations of $f$ has a subsequence that converges uniformly. Moreover (\cite{gottschalk}, 4.48) implies that the following conditions are equivalent:
\begin{enumerate}
\item $f \in B(\mathbb R^n).$
\item $\overline {O(f)}$ is a compact minimal closed subset of $C_b(\mathbb R^n)$ invariant under $\pi^x, x \in \mathbb R^n.$
\item there is a unique group structure on 
$\overline {O(f)}$ which makes it a compact topological group  and the map $\psi : \mathbb R^n \mapsto G$ defined by
$\psi(x) := \pi^xf$ has a dense image.
\end{enumerate}
If $f \in B(\mathbb R^n)$ define  $G := \overline {O(f)}$ and 
$F \in C(G)$ by $F(h) := h(0)$ for $h \in G.$ Then $f = F \circ \psi.$ The Stone-Weierstrass theorem (\cite{rudin1}, A14) implies that the algebra of trigonometric polynomials $\mathcal T(G)$ on $G$ is dense in $C(G)$ so $F$ can be uniformly approximated by $p \in \mathcal T(G)$ hence $f$ can be uniformly approximated by $p \circ \psi \in \mathcal T_n$ which implies Bohr's fundamental theorem. These results illustrate the utility of topological dynamics
\end{remark}
If $f \in B(\mathbb R^n)$ and
\begin{equation}\label{summable1}
\sum_{\omega \in \Omega(f)} |\mathcal F_B(f)(\omega)| < \infty
\end{equation}
then 
\begin{equation}\label{FB expansion}
f = \sum_{\omega \in \Omega(f)} \mathcal F_B(f)(\omega) \, \widehat \omega,
\end{equation}
where $\widehat \omega$ is defined in (\ref{hatw1}). The right side of (\ref{FB expansion}) converges uniformly independent of the ordering of terms.
(\cite{bohr2}, p. 89, Convergence Theorem) shows that this always occurs whenever 
$\Omega(f)$ is linearly independent over $\mathbb Q.$
If (\ref{summable1}) fails then 
the right side of (\ref{FB expansion}) must be computed as a limit of Fejer sums as explained in (\cite{bohr2}, p. 66). 
The convolution of $f \in B(\mathbb R^n)$ and $\xi \in C_c(\mathbb R^n),$ \\
$(\xi*f)(x) := \int_{y \in \mathbb R^n} f(x-y)\xi(y)dy$
is in $B(\mathbb R^n)$ and
\begin{equation}\label{FBconvf}
	\mathcal F_B(\xi*f) = \widehat \xi\, \mathcal F_B(f).
\end{equation}
\begin{defi}\label{BAPmu}
A Radon measure $\mu$ is a Bohr (resp. Besicovitch) almost periodic measure if for every $\xi \in C_c(\mathbb R^n),$ the convolution 
$\xi*\mu := \, <\mu, \xi(x - \cdot)>$ is a Bohr (resp. Besicovitch) almost periodic function. 
A Bohr (resp. Besicovitch) almost periodic distribution is a tempered Radon measure $\mu$ such that for every $\xi \in S(\mathbb R^n),$ $\xi*\mu$ is a Bohr (resp. Besicovitch) almost periodic function.
\end{defi}
If $\mu$ is a Bohr almost periodic measure define 
$\mathcal F_B(\mu) : \mathbb R^n \mapsto \mathbb C$ by
\begin{equation}\label{FBmu}
\mathcal F_B(\mu)(\omega) := \frac{\mathcal F_B(\xi*\mu)(\omega)}{\widehat \xi(\omega)}
\end{equation}
where $\xi \in C_c(\mathbb R^n)$ and $\widehat \xi(\omega) \neq 0.$ (\ref{FBconvf}) ensures that $\mathcal F_B(\mu)$ is well defined. The spectrum 
\begin{equation}\label{specmu}
\Omega(\mu) := \{\omega \in \mathbb R^n : \mathcal F_B(\mu)(\omega) \neq 0\}.
\end{equation}
is countable since there exist $\xi \in C_c(\mathbb R^n)$ such that $\widehat \xi(\omega) \neq 0$ for all $\omega.$
Meyer (\cite{meyer4}, 2nd sentence) calls Radon measures that are are countable linear combinations point measures purely atomic measures. Meyer \cite{meyer8} and Favorov \cite{favorov8} define a tempered distribution a $\mu$ to be a Poisson measure if both $\mu$ and $\widehat \mu$ are purely atomic measures. 
\begin{lem}\label{meyer} 
If $\mu$ is a Bohr almost periodic measure then
$\widehat \mu \in \mathcal R(\mathbb R^n)$ iff
\begin{equation}\label{summable}
\sum_{\omega \in \Omega(\mu) \cap B_n(0,r)} 
|\mathcal F_B(\mu)(\omega)| < \infty, \ \ r > 0.
\end{equation}
Then 
\begin{equation}
\widehat \mu = \sum_{\omega \in \Omega(\mu)} 
\mathcal F_B(\mu)(\omega)\, \delta_\omega
\end{equation}
is a purely atomic measure. Therefore, if $\mu$ is a Bohr almost periodic purely atomic measure, then $\mu$ is a Poisson measure iff (\ref{summable}) holds.
\end{lem}
Proof. Meyer proved this (\cite{meyer5}, Theorem 5.5).
%
%
\begin{remark}\label{favorov} 
Favorov (\cite{favorov4}, Lemma 1) proved that if $\mu$ is a tempered measure, $\widehat \mu$ is a pure 
point measure, and 
$|\widehat f|(B_n(0,r)) = O(r^\gamma)$ for some $\gamma > 0,$ then  $\mu$ is an almost periodic distribution and for every 
$\psi \in \mathcal S(\mathbb R^n),$ $\mathcal F_B(\psi*\mu)$ converges absolutely.
\end{remark}
We record the following observation without proof. 
\begin{remark}\label{BAPdist} 
If $\mu$ is a Bohr almost periodic measure on $\mathbb R$ and (\ref{summable}) fails, then there exists $h \in C(\mathbb R)$ such that 
$\widehat \mu = h^{\, \prime \prime}$ and $h^{\, \prime}$ is differentiable and has derivative $ = 0$ except on $\Omega(\mu)$
where $h^{\, \prime}$ has a jump equal to $\mathcal F_B(\mu)(\omega)$ at every $\omega \in \Omega(\mu).$ In higher dimensions $\widehat \mu$ is more complicated.
\end{remark}
\begin{defi}\label{tbmu}
A Radon measure $\mu$ on $\mathbb R^n$ is translation bounded if
there exists $M > 0$ such that
\begin{equation}\label{transbdd}
\sup_{x \in \mathbb R^n} |\mu|(B_n(x,1)) \leq M.
\end{equation}
\end{defi}
Since $B_n(0,r)$ is contained in a union of $O(r^n)$ balls of the form $B_n(x,1),$ Lemma \ref{trm2} implies that every translation bounded 
Radon measure is a tempered Radon measure. 
%
%
\begin{prop}\label{prop1}
If $\mu \in \mathcal R(\mathbb R^n)$ then $\mu$ is a Bohr almost periodic measure iff $\mu$ is translation bounded and is a Bohr almost periodic distribution.
\end{prop}
Proof. Let $\mu$ be a Bohr almost periodic measure. Let $V$ be the Banach space of consisting of $\xi \in C_c(\mathbb R^n)$ whose supports are contained in 
$\overline {B_n(0,1)}$ and let $V_1 := \{\xi \in V : ||\xi||_\infty \leq 1\}.$
For every $p \in \mathbb R^n$ define $L_p \in V^*$ by
\begin{equation}
	L_p(\xi) := (\xi*\mu)(p).
\end{equation}
Clearly each $L_p$ is a bounded lineear functional on $V$ since
\begin{equation}
	|L_p(\xi)| \leq |\mu|\, (B_n(p,1)) \times ||\xi||_\infty.
\end{equation}
Furthermore, since $\xi*\mu \in B(\mathbb R^n) \subset C_b(\mathbb R^n),$
\begin{equation}
	\sup_{p \in \mathbb R^n} |L_p(\xi)| =
	\sup_{p \in \mathbb R^n} |(\xi*\mu)(p)| < \infty, \ \ \xi \in V.
\end{equation}
Therefore the Banach-Steinhaus theorem (\cite{rudin2}, 5.8)
implies that there exists $M > 0$ such that 
\begin{equation}
	||L_p|| :=
	\sup_{\xi\in V_1} |L_p(\xi)| = |\mu|(B_n(p,1)) < M, \ \ p \in \mathbb R^n,
\end{equation}
hence $\mu$ is translation bounded. Let $\xi \in \mathcal S(\mathbb R^n).$ We conclude the proof of the only if part by proving that $\xi*\mu \in B(\mathbb R^n).$ Let $\xi_j \in \mathcal S(\mathbb R^n)$ be a sequence such that $\xi - \xi_j$ equals $0$ on $B_n(0,j)$ and $|\xi(x)-\xi_j(x)| < C(1+r)^{-n-2}$ for $|x| > r.$ Then $\xi*\mu \in C_b(\mathbb R^n)$ and $\xi_j*\mu$ converges uniformly to $\xi*\mu.$ Since each $\xi_j*\mu \in B(\mathbb R^n)$ and $B(\mathbb R^n)$ is closed under uniform limits, $\xi*\mu \in B(\mathbb R^n).$
To prove the if part assume that $\mu$ is a translation bounded measure and a Bohr almost periodic distribution. Let $f \in C_c(\mathbb R^n).$ We conclude the if part of the proof by showing that
$f*\xi \in B(\mathbb R^n).$ Let $\xi_j \in \mathcal S_c(\mathbb R^n)$ such that $||f - \xi_j||_1 \rightarrow 0.$ Then $f*\mu$ is a bounded continuous function and is the uniform limit of $\xi_j*\mu \in B(\mathbb R^n)$ hence $f*\mu \in B(\mathbb R^n).$
\\ \\
Sergey Favorov proved (\cite{favorov9}, Theorem 1) that every non-negative Poisson measure is a Bohr almost periodic measure.  Proposition \ref{prop1} implies that his result gives the following result as a special case. We provide a proof to make the paper more self-contained.
is
%
%
\begin{prop}\label{prop2}
Every positive FQ on $\mathbb R^n$ 
is translation bounded and is a Bohr almost periodic distribution hence is a Bohr almost periodic measure.
\end{prop}
Proof. Let $\mu$ be a positive FQ so $\widehat \mu = \sum_{s \in S} a_s\, \delta_s,$ where $S$ is discrete, $a_s \in \mathbb C^*,$ and 
$|\, \widehat \mu \, | = \sum_{s \in S} |a_s|\, \delta_s$ 
is a tempered distribution. Therefore there exists  $\beta > 0$ such that
\begin{equation}\label{eqq1}
F(r) := |\, \widehat \mu \, |(B_n(0,r)) = \sum_{|s| < r} |a_s| = O(r^\beta), \ \ r \rightarrow \infty.
\end{equation}
Let $\xi \in \mathcal S(\mathbb R^n).$ Then 
$\xi(x - \cdot) = \widehat h$ where
$h(y) = e^{2\pi i x \cdot y} \widehat \xi(y),$ hence
\begin{equation}\label{eqq2}
(\mu*\xi)(x) = \, <\mu, \widehat h> \, = \, <\widehat \mu, h> \, 
= \sum_{s \in S} a_s \, e^{2\pi i s \cdot x} \widehat \xi(s).
\end{equation}
Since $\widehat \xi \in \mathcal S(\mathbb R^n),$ there exists $C > 0$ such that
$|\widehat \xi(s)| < C |s|^{-\beta - 1}$ for $|s| \geq 1.$ Therefore there exist
$C_0, C_1 > 0$ with
\begin{equation}\label{eqq3}
\sum_{s \in S} |a_s| \, |\widehat \xi(s)| \leq C_0 + C_1 \int_1^\infty r^{-\beta - 1} dF(r) < \infty.
\end{equation}
This proves $\mu*\xi \in B(\mathbb R^n) $ and has an absolutely converging Fourier-Bohr series so $\mu$ is a Bohr almost periodic distribution. Choose $\xi \in S_c(\mathbb R^n$ values in $[0,1]$ whose support is in $B_n(0,2)$ and equals $1$ on $B_n(0,1).$ Since $\mu*\xi \in B(\mathbb R^n)$ it is bounded and $\mu$ positive implies $|\mu| = \mu$ hence
\begin{equation}\label{eqq4}
\sup_{x \in \mathbb R^n} |\mu|(B_n(x,1)) \leq \sup_{x \in \mathbb R^n} (\xi*\mu)(x) < \infty
\end{equation}
so $\mu$ is translation bounded. Proposition \ref{prop1} then implies that $\mu$ is a Bohr almost periodic measure.
\subsection{Bohr and Besicovitch Almost Periodic Multisets}\label{subsec2.5}
Krein and Levin \cite{kreinlevin} introduced almost periodic multisets in $\mathbb C$ to study zeros of holomorphic almost periodic functions. Favorov and Kolbasina \cite{favorovkolbasina} studied 
them in $\mathbb R^n.$  In this subsection we discuss them in $\mathbb C^n$ to study zeros of holomorphic maps in Subsection \ref{subsec2.6}.
\begin{defi}
A multiset in $\mathbb C^n$ is a pair 
$\alpha = (\Lambda_\alpha,m_\alpha)$ where 
$\Lambda_\alpha \subset \mathbb C^n$ is discrete and 
$m_\alpha : \mathbb C^n \mapsto \mathbb Z_{\geq 0}$ has support 
$\Lambda_\alpha.$ We let $\mathcal S(\mathbb C^n)$ denote the set
of multisets in $\mathbb C^n$ and $\mathcal S_1(\mathbb C^n)$
denote its subset where $m_\alpha = 1$ on $\Lambda_\alpha.$ If 
$\Lambda_\alpha \subset \mathbb R^n$ we say $\alpha$ is a multiset 
in $\mathbb R^n.$
\end{defi}
$\alpha \in \mathcal S(\mathbb C^n)$ can be represented by a sequence $\lambda_k \in \Lambda_\alpha, k \in \mathbb Z_+$ where for every $\lambda \in \Lambda_\alpha,$ $\lambda_k = \lambda$ for $m_\alpha(\lambda)$ values of $k.$ Two sequences $\lambda_k$ and $\mu_k$ represent the same multiset iff there exists a permutation
$\sigma : \mathbb Z_+ \mapsto \mathbb Z_+$ such that
$\lambda_k = \mu_{\sigma(k_)}.$ This defines equivalence classes of $\mathbb C^n$--valued sequences.  
\\ \\
Define the distance between multisets
\begin{equation}\label{distms}
	d(\alpha,\beta) := \inf_\sigma \sup_{k \in \mathbb Z_+} |\lambda_k - \mu_{\sigma(k)}|,
\end{equation}
where $\alpha$ is represented by a sequence $\lambda_k,$ $\beta$ is represented by a sequence $\mu_k,$ and $\sigma$ ranges over all permutations of $\mathbb Z_+.$ We observe that the metric space $(\mathcal S(\mathbb C^n),d)$ is complete since every Cauchy sequence converges. 
\begin{remark}\label{Hausdistance}
The restriction of $d$ to $\mathcal S_1(\mathbb C^n)$ is the
Hausdorff distance explained in (\cite{munkres}, p. 279, Exercise 7)
\end{remark}
For $x \in \mathbb R^n$ define the translation $\pi^x : \mathcal S(\mathbb C^n) \mapsto \mathcal S(\mathbb C^n)$
by 
$$\pi^x(\Lambda_\alpha, m_\alpha) := 
(\Lambda_\alpha + x, m_\alpha(\, \cdot - x)).$$
Since translations are isometries, the dynamical system is uniformly equicontinuous. 
Following the logic used to analyze Bohr almost periodic functions in Subsection \ref{subsec2.4},
we define $\alpha$ to be Bohr almost periodic if for every $\epsilon >0$ the set of
$\epsilon$-almost periods 
$$S_\alpha(\epsilon) := \{x \in \mathbb R^n : d(\alpha,\pi^x(\alpha)) \leq \epsilon\}$$ 
is syndetic. Therefore $\alpha$ is Bohr almost periodic iff its orbit closure $\overline {O(\alpha)}$ is compact and then 
$\overline {O(\alpha)}$ admits a unique compact group structure such that the map $\psi : \mathbb R^n \mapsto \overline {O(\alpha)},$
defined by $\psi(x) := \pi^x(\alpha),$ is a continuous homomorphism
with a dense image. Thus the pair $(\overline {O(\alpha)}, \psi)$ is a compactification of $\mathbb R^n$ in Definition \ref{compactification}.
Define the sum of multisets by
$$(\Lambda_\alpha, m_\alpha) \oplus (\Lambda_\beta, m_\beta)
:= (\Lambda_\alpha \cup \Lambda_\beta,m_\alpha+m_\beta).$$
Associate to a multiplicity set
$\alpha = (\Lambda_\alpha, m_\alpha)$ the unique measure 
\begin{equation}\label{mualpha}
	\mu_\alpha := \sum_{\lambda \in \Lambda_\alpha} m_\alpha(\lambda) \, \delta_\lambda.
\end{equation}
Kolbasina and Favorov \cite{favorovkolbasina, kolbasina}  proved that a multiset $\alpha$ in $\mathbb R^n$ is Bohr almost periodic iff
$\mu_\alpha$ is a Bohr almost periodic measure and then its density
\begin{equation}\label{density1}
	\Delta(\Lambda,m_\Lambda) := 
\lim_{r \rightarrow \infty} \frac{r^{-n}}{c_n} \, \mu_\alpha(B_n(y,r)), \ \ y \in \mathbb  R^n
\end{equation}
exists, is independent of $y,$ and the sum converges uniform in $y.$ 
We define a multiset $\alpha$ in $\mathbb R^n$ to be Besicovitch almost periodic if $\mu_\alpha$ is a Besicovitch almost periodic measure.
Favorov\cite{favorov3}  characterized Besicovitch almost periodic multisets in $\mathbb R^n$ by their translation properties.
\begin{defi} 
For
$\alpha := (\Lambda_\alpha, m_\alpha) \in \mathcal S(\mathbb C^n)$
define
$$b(\alpha) : = \inf_{\lambda_1, \lambda_2 \in \Lambda_\alpha, \lambda_1 \neq \lambda_2} |\lambda_1-\lambda_2| > 0.$$
$\alpha$ is uniformly discrete if 
$\alpha \in \mathcal S_1(\mathbb C^n)$ and $b(\alpha) > 0.$
$\alpha$ is relatively uniformly discrete if it 
is the sum of a finite number of uniformly discrete multisets. 
The index of $\alpha$ is
\begin{equation}\label{ind}
	ind(\alpha) := 
		\inf_{\epsilon > 0} \sup_{z \in \mathbb C^n} 
		\sum_{\lambda \in \Lambda_\alpha \cap B_{c,n}(z,\epsilon)} m_\alpha(\lambda).
\end{equation}
\end{defi}
We observe that $\alpha$ is relatively uniformly discrete iff $ind(\alpha) < \infty$ and then $\alpha$ equals the sum of $ind(\alpha)$ uniformly discrete multisets. Furthermore, all Bohr almost periodic multisets in $\mathbb R^n$ are relatively uniformly discrete.
\begin{remark}\label{covering} 
Assume that $\alpha \in \mathcal S_1(\mathbb R^n)$ is uniformly discrete, Bohr almost periodic, and $\dim_{\mathbb Z} \Omega(\mu) = m < \infty$ 
and let $(\overline {O(\alpha)}, \psi)$ be the associated compactification of $\mathbb R^n.$ Then (\cite{lawton}, Theorem 5) implies that 
$\overline {O(\alpha)} \cong\mathbb T^m$ 
and $\overline {\psi(\Lambda_\alpha)},$ called the compactification of $\Lambda_\alpha,$ is homeomorphic to a finite union of pairwise disjoint submanifolds of $\mathbb T^m$ each homeomorphic
to $\mathbb T^{m - n}$ and embdede transverse to the foliation induced by $\psi$ and inducing an injective homotopy. 
\end{remark}
\subsection{Holomorphic Maps and Multiplicity of their Roots}\label{subsec2.7}
A domain $D \subset \mathbb C^n$ is a nonempty open connected set.
A map $F = [f_1,...,f_n]^T: D \rightarrow \mathbb C^n$ is holomorphic if each component $f_j$ is holomorphic. 
A root $\lambda \in \Lambda(F)$ of $F$ is isolated if there exists 
$\epsilon > 0$ such that 
$\overline {B_{c,n}(\lambda,\epsilon)} \subset D$ and
$\Lambda(F) \cap \overline {B_{c,n}(\lambda,\epsilon)} = \{\lambda\}.$ Assume that $\lambda$ is an isolated root of $F.$ For $n = 1,$ $F = f$ and the multiplicity of $\lambda$ equals its logarithmic integral
\begin{equation}
m_\lambda(f) = \frac{1}{2\pi i} \int_{C_\lambda} \frac{df}{f}
\end{equation}
where $C_\lambda$ is the one-cycle consisting of a circle with center at $\lambda$ bounding a closed disc contained in $D$ and containing only the root $\lambda$ 
and oriented so $d(arg \, z)$ is positive (counterclockwise orientation). For $n > 1$ the multiplicity can be computed by a Bochner-Martinelli integral over a sufficiently small $(2n-1)$-dimensional sphere centered at $\lambda$ (\cite{tsikh}, p. 16, Eqn. 1), (\cite{shabat}, p. 286, Eqn. 3). 
Alternatively, its multiplicity can be computed by an integral 
\begin{equation}\label{growint}
m_\lambda(F) = \left(\frac{1}{2\pi i}\right)^n \int_{\Gamma_\lambda} 
\frac{df_1}{f_1} \wedge \cdots \wedge \frac{df_n}{f_n},
\end{equation}
over the $n$-cycle
\begin{equation}\label{Gcycle}
\Gamma_\lambda := \{z \in \mathcal U_\lambda : |f_j(z)| = \epsilon_j, \, j = 1,...,n, \},
\end{equation}
where $\mathcal U_\lambda$ is a neighborhood of $\lambda$ with
$\mathcal U_\lambda \cap \Lambda(P) = \{\lambda\},$
$\epsilon_j > 0$ are sufficiently small to ensure that $\Gamma_\lambda \subset \mathcal U_\lambda,$ 
and $\Gamma_\lambda$ is oriented so that the differential form
$d(arg f_1) \wedge \cdots \wedge d(arg f_n) \geq 0.$ 
We call $\Gamma_\lambda$ a Grothendieck cycle since (\ref{growint}) is the Grothendieck residue of the Jacobian $J(f,z)$ (\cite{tsikh}, p. 16). The Bochner-Martinelli integral representation of multiplicity gives:
\begin{lem}\label{rouche} 
(Rouche) If $D \subset \mathbb C^n$ is a bounded domain with Jordan boundary $\partial D,$ and if $f, g : \overline D \rightarrow \mathbb C^n$ are continuous, holomorphic on $D,$ and
$|g| < |f|$ on $\partial D,$ then $f+g$ has as many zeros in $D$ as 
$f.$
\end{lem}
Proof. Follows from (\cite{shabat}, Theorem 2, p. 288).
\begin{defi}\label{tubedomain}
A (horizontal) tube domain is a subset $T = \mathbb R^n + iK \subset \mathbb C^n$ where $K \subset \mathbb R^n$ is nonempty, open, connected and bounded. 
\end{defi}
The set of bounded holomorphic functions on $T$ is a Banach space with norm $||F||_T := \sup_{z \in T} |F(z)|.$  For $\Lambda(F)$ discrete define the multiset $div\, F := \alpha$ where 
$\Lambda_\alpha := \Lambda(F)$ and 
$m_\alpha : T \mapsto \mathbb Z_{\geq 0}$ is supported on $\Lambda(F)$ and $m_\alpha(\lambda) := m_\lambda(F)$ for 
$\lambda \in \Lambda(F).$ We also identifty $div\, F$ with the Radon measure $\mu_\alpha.$
\begin{defi}\label{HAP}
A holomorphic map $F : T \mapsto \mathbb C^n$ is Bohr almost periodic if it is bounded and for every $\epsilon > 0$ there exists a syndetic $S_F(\epsilon) \subset \mathbb R^n$ such that 
\begin{equation}
	||F(\cdot - \tau)-F||_T \leq \epsilon, \ \ \tau \in S_F(\epsilon).
\end{equation}
$B(T)$ denotes the Banach space of holomorphic Bohr almost periodic maps $F :T \mapsto \mathbb C^n.$
\end{defi}
Roots are counted with multiplicities. 
If $F \in B(T)$ and $\Lambda(F)$ is discrete then
$div\, F$ is a multiset and a Radon measure and we define
the lower density of $div\, F$ by
\begin{equation}\label{dendivF}
	\Delta_{low}(div\, F) := \liminf_{r \rightarrow \infty} \frac{r^{-n}}{c_n} \sum_{\lambda \in \Lambda(F) \cap B_{c,n}(0,r)}  m_\lambda(F).
\end{equation}
\begin{cor}\label{Rouchcor1} 
$\Delta_{low}(div\, F) > 0.$
\end{cor}
Proof.  Let $\lambda \in \Lambda(F),$ $r > 0$ with 
$\overline {B_{c,n}(\lambda,r)} \subset T,$  
$\partial B_{c,n}(\lambda,r) \cap \Lambda(F) = \emptyset.$ Since
$\epsilon := \frac{1}{2} \min_{z \in \partial B_{c,n}(\lambda,r)} |F(z)| > 0$ 
there exists a syndetic 
$S_\epsilon(F) \subset \mathbb R^n$ with
$||F(\cdot - \tau) - F||_T \leq \epsilon, \ \ \tau \in S_\epsilon(F).$
$F$ has $k := div\, F(B_{c,n}(\lambda,r)) > 0$ roots in 
$B_{c,n}(\lambda,r)$ and 
$|F(z - \tau) - F(z)| < |F(z)|, z \in \partial B_{c,n}(\lambda,r)$ 
so Lemma \ref{rouche} implies that $F(\cdot - \tau)$ has $k$ roots in
$B_{c,n}(\lambda,r)$ and $F$ has $k$ roots in 
$B_{c.n}(\lambda+\tau,r).$
Since $S_\epsilon(F)$
is syndetic it must contain a subset $S$ with lower density
$\Delta_{low}(S) > 0$
such that the balls $B_{c,n}(\lambda - s,r), s \in S$ are pairwise disjoint. This implies that
$\Delta_{low}(div\, F) \geq k \Delta_{low}(S) > 0$
and concludes the proof.
\\ \\
For $\epsilon > 0$ define the compact set $K_\epsilon := \{x \in K : d(x, \mathbb R^n \backslash K) \geq \epsilon\}$ and the closed set 
$T_\epsilon := \{z \in T : d(z,\Lambda(F)) \geq \epsilon\} \cap (\mathbb R^n + iK_\epsilon).$ 
\begin{cor}\label{Rouchcor2} 
If $F \in B(T)$ and $\Lambda(H)$ is discrete for every $H \in \overline {O(F)},$ then for every $\epsilon > 0,$
$\inf_{z \in T_\epsilon} |F(z)| > 0.$ 
\end{cor}
Proof.  Assume to the contrary that there exists $\epsilon > -$ with $\inf_{z \in T_\epsilon} |F(z)| = 0.$ 
Then there exists a sequence
$z_j = x_j + iy_j \in T_\epsilon$ with $F(z_j) \rightarrow 0.$ 
Since both  
$\overline {O(F)}$ and $K_\epsilon$ are compact, we may replace $z_j$ with a subsequence satisfying the following property: there exist $y_0 \in K_\epsilon$ and $H \in \overline {O(F)}$ such that
$$y_j \rightarrow y_0 \hbox{ and } \ \ ||F(\cdot + x_j)-H||_T. \rightarrow 0.$$
Therefore
$$H(iy_0) = \lim_{j \rightarrow \infty} H(iy_j) 
= \lim_{j \rightarrow \infty} F(iy_j + x_j)
= \lim_{j \rightarrow \infty} F(z_j) = 0.$$
Since $\Lambda(H)$ is discrete there exists $0 < r < \epsilon/2$ 
such that
$$\overline  {B_{c,n}(iy_0,r)} \subset T \hbox{ and } 
\partial B_{c,n}(iy_0,r) \cap \Lambda(H) = \emptyset.$$
Lemma \ref{rouche} implies that for sufficiently large $j,$
$F(\cdot + x_j)$ has a root $\lambda_j \in B_{c,n}(iy_0,r).$
Since $z_ = x_j+iy_j \in T_\epsilon$ and $x_j+\lambda_j \in \Lambda(F),$
\begin{equation}\label{lowbdd}
|\lambda_j - y_j| \geq \epsilon.
\end{equation}
We may also choose $j$ so that $|y_j - y_0| < r$ hence
since $\lambda_j + x_j \in B_{c,n}(iy_0,r),$
$$|\lambda_j - iy_j| \leq |\lambda_j - iy_0| + |y_0 - y_j| < 2r  < \epsilon.$$
This contradicts (\ref{lowbdd}) and concludes the proof.
\begin{cor}\label{Rouchcor3} 
If $F \in B(T),$ $\Lambda(F) \neq \emptyset,$ and for every $H \in \overline {O(F)}$ the set $\Lambda(H)$ is discrete, then 
\begin{enumerate}
\item If $C \subset T$ is compact, then $\sup_{x \in \mathbb R^n} |\Lambda(F) \cap (C+x)| < \infty.$
\item $div\, F$ is a Bohr almost periodic multiset.
\item $\Omega(div\, F)$ is contained in the group generated by 
$\Omega(F).$
\end{enumerate}
\end{cor}
Proof. To prove 1. assume to the contrary that there exists a sequence $x_j \in \mathbb R^n$ such that 
$$|\Lambda(F(\cdot + x_j)) \cap C| = |\Lambda(F )\cap  (C+x_j)| \rightarrow \infty.$$  Since $F \in B(T),$ 
$\overline {O(F)}$ is compact so $x_j$ can be replaced by a 
subequence such that $||F(\cdot + x_j) - H||_T \rightarrow 0.$ 
Since $F$ is bounded the partial derivatives of its components 
are bounded hence $H$ has an infinite number of zeros in $C.$ 
This contradicst the assumption that $\Lambda(H)$ is discrete and 
concludes the proof.
\\
To prove 2. we first observe that 1. implies there exists $r > 0$ 
and $N \in \mathbb Z_+$ such that each connected component of 
$T_{r/2} \cap \cup_{\lambda \in \Lambda(F)} \overline {B_{c,n}(\lambda,r)}$ 
is the union of at most $N$ sets
$T_{r /2} \cap \overline {B_{c,n}(\lambda,r)}.$ 
Let $D_k, k \in \mathbb Z_+$ be interiors of these connected components. The diameter of each $D_k$ is $\leq 2rN.$ Corollary \ref{Rouchcor2} gives
$$\delta := \inf_{k \in \mathbb Z_+, \, z \in \partial T_{r/2} \cap D_k} |F(z)|  > 0 .$$
There exists a syndetic $S_F(\delta)$ with
$$||F(\cdot+\tau) - F||_T < \delta, \ \ \tau \in S_F(\delta).$$ 
Therefore Lemma \ref{rouche} implies that 
$$div\, F(D_k + \tau) = div\, F(\cdot + \tau)(D_k = 
(\tau + div\, F), \  \tau \in S_F(\delta).$$ 
Since the diameter of each $D_k$ is $\leq 2rN$ this implies
$d(div\, F +\tau ,div\, F) \leq Nr$ 
hence $div\, F$ is a Bohr almost periodic multiset. 
\\ \\
To prove 3. define $\alpha := div\, F$ and observe that 2. implies that 
the pair $(\overline {O(\alpha)},\psi_2),$ where $\psi_2 : \mathbb R^n \mapsto O(\alpha)$ is defined by $\psi_2(x) := \pi^x(\alpha),$ is a compactification of $\mathbb R^n$ as is
Definition \ref{compactification}. Also the pair
$(\overline {O(F)},\psi_1),$ where $\psi_1 : \mathbb R^n \mapsto O(F)$ is defined by $\psi_1(x) := F(\cdot + x),$ is a compactification of $\mathbb R^n.$ Lemma \ref{fF} implies that the spectrum
$\Omega(F)$ generates the spectrum $\Omega(\overline {O(F)},\psi_1)$ of the comactification and hence is contained in it. A similar argument gives that $\Omega(\alpha)$ generates the spectrum $\Omega(\overline {O(\alpha)},\psi_2)$ of the compactification and hence is contained in it.
Lemma \ref{rouche} implies that the homomorphism $\eta_1 : O(F) \mapsto O(\alpha)$ defined by $\eta_1(F(\cdot + x) := \pi^x(\alpha)$ 
is uniformly continuous with respect to the metric 
$|| \cdot ||+T$ on $O(F)$ and the metrix $d$ on 
$O(\alpha).$ Therefore $\eta_1$ extends to a 
continuous homomorphism 
$\eta : \overline {O(F)} \mapsto \overline {O(\alpha)}$ and $\eta$ 
is an equivariant map between minimal dynamical systems hence 
$\overline {O(F)} \cong \overline {O(\alpha)}$ and 
$\Omega(\overline {O(F)},\psi_1) = \Omega(\overline {O(\alpha)},\psi_2).$ 
Therefore 
$$\Omega(\alpha) \subset \Omega(\overline {O(\alpha)},\psi_2) = 
\Omega(\overline {O(F)},\psi_1) = \text{ group generated by } \Omega(F).$$
We note that the proof of Lemma \ref{fF} does not require any result from Corollary \ref{Rouchcor3} so the proof of 3. does not employ circular
reasoning. The following example shows that the inclusion in 3. can be strict. If $T \subset \mathbb C$ is a tube domain containing $\mathbb R$ then $\sin 2\pi z \in B(T),$ $\Omega(\sin 2\pi z) = \{-1.1\}$ and Poisson's summation formula implies $\Omega(div\, sin 2\pi z) = \mathbb Z.$
\begin{remark}\label{favorov} 
The study of almost periodic divisors and zero sets of almost periodic holomorphic functions was started by Lev Isaakovich Ronkin \cite{ronkin} and his co-authors Sergey Yuryevich Favorov and Alexander Ju. Rashkovskii \cite{favorovRR1, favorovRR2, favorovRR3}.
\end{remark}
\subsection{Representation of Trigonometric Maps}\label{subsec2.6}
For $Q_1, Q_2 \in \mathcal L_m^n$ we write $Q_1 \sim Q_2$ if the exists $A \in GL(m,\mathbb Z)$ such that $Q_2 = Q_1 \circ \widetilde A.$ This is an equivalence relation. 
\begin{defi}
$Q \in \mathcal L_m^n$ is minimal if $\mathbb Z(\Omega(Q)) = \mathbb Z^m.$
\end{defi}
\begin{prop}\label{prop3} 
If $Q \in \mathcal L_m^n$ and $ M \in \mathbb R^{m \times n},$ then
\begin{equation}\label{QM}
P := Q \circ \rho_m \circ M \in \mathcal T_n^n,
\end{equation}
where $\rho_m$ is defined in (\ref{rhom}), and
\begin{equation}\label{OmegaQM}
\Omega(P) = M^T \, \Omega(Q).
\end{equation}
If $P \in \mathcal T_n^n$ and $m = \dim_{\mathbb Z} \Omega(P)$ then there exists a minimal $Q \in \mathcal L_m^n$ and $M \in \mathbb R^{m \times n}$ such that $P = Q \circ \rho_m \circ M.$
For every minimal $Q_1 \in \mathcal L_m^n$ and $M_1 \in \mathbb R^{m \times n}$ that satisfy $P = Q_1 \circ \rho_m \circ M_1$ there exists $A \in GL(m,\mathbb Z)$ such that $M_1 = A^{-T}M$ and $Q_1 = Q \circ \widetilde A.$
\end{prop}
Proof. The first assertion follows from a direct computation. 
Let $m := \hbox{dim}_{\mathbb Z}(\Omega(P))$ and $W \in \mathbb R^{n \times m}$ whose columns form a $\mathbb Z$-basis for $Z(\Omega(P))$ and $M := W^T.$ 
For every $\omega \in \Omega(P)$ there exists a unique $\ell_\omega \in \mathbb Z^m$ such that $\omega = W\ell_{\omega}.$ 
Then for every entry $p_j \in \mathcal T_m$ of $P,$ (\ref{hatw2}) implies that
$$p_j = \sum_{\omega \in \Omega(p)} c_\omega \widehat {\omega} 
= q_j \circ \rho_m \circ M$$
where $q_j \in \mathcal L_m$ is defined by
$$q_j = \sum_{\omega \in \Omega(p_j)} c_\omega \widetilde {\ell_\omega}$$
and $\widetilde {\ell_\omega}$ is defined following (\ref{tildeN}).
Therefore $P = Q \circ \rho_m \circ M$ where $Q := [q_1,...,a_m]^T.$
Since $\Omega(P) = M^T \Omega(Q),$ 
$W\mathbb Z^m = \mathbb Z(\Omega(P)) = M^T \mathbb Z(\Omega(Q)) = W \mathbb Z(\Omega(Q)).$ 
Since the column of $W$ are $\mathbb Z$-independent, 
$\mathbb Z(\Omega(Q)) = \mathbb Z^m$ so $Q$ is minimal. To prove the third assertion note that
$W\mathbb Z^m = \mathbb Z(M_1^T \Omega(Q_1)) = M_1^T \mathbb Z^m$ hence there exist
$A, B \in \mathbb Z^{m \times m}$ with $W = M_1^T B$ and $M_1^T = WA$ and 
Then $W(I_m - AB) = 0.$ Since the columns of $W$ are $\mathbb Z$-independent, $AB = I_m$ hence
$A, B \in GL(m,\mathbb Z)$ hence $M_1 = A^{-T}M.$ Let $\psi := \rho_m \circ M.$ Then (\ref{rhomtildeN}) implies
$$P = Q \circ \psi = Q_1 \circ \rho_m \circ M_1 = Q_1 \circ \rho_m \circ  B^TM = Q_1 \circ \widetilde B \circ \psi.$$
Lemma 1 implies that $\psi$ has a dense image hence $Q_1 \circ \widetilde B = Q$ so $Q_1 = Q \circ \widetilde A.$
\begin{cor}\label{corProp3} corollary4
If $P \in \mathcal T_n^n$ and $T$ is a tube domain then $P|_T \in B(T).$
\end{cor}
Proof. For $\theta \in [-\pi,\pi]$ and $s \in \mathbb R,$
$|\theta| \leq \frac{\pi}{2} \, |e^{i\theta} - 1| \leq \frac{\pi}{2} |\theta|,$ hence
\begin{equation}\label{stheta}
|e^{is\theta} - 1| \leq \frac{|s|\pi}{2} \, |e^{i\theta} - 1|.
\end{equation}
Proposition \ref{prop3} implies there exists $m \geq n,$ $Q \in \mathcal L_m^n$ and
$M \in \mathbb R^{m \times n}$ such that $P = Q \circ \rho_m \circ M.$ 
Let $p$ is a component of $P$ and $q$ the corresponding component of $Q.$ 
Then $p = q \circ \rho_m \circ M.$ 
Let $q = \sum_{\ell \in \Omega(q)} c_\ell \widetilde \ell.$
For $z \in T,$ $\tau \in \mathbb R^n,$ the triangle inequality gives
\begin{equation}\label{pbound1}
|p(z + \tau) - p(z)| \leq C_q \max_{\ell \in \Omega(q)} |e^{2\pi i \ell^TM\tau}-1| 
\end{equation}
where
\begin{equation}\label{Cq}
C_q := \sup_{y \in K} \sum_{\ell \in \Omega(q)}  |c_\ell| \, e^{-2\pi \ell^T My}
\end{equation}
(\ref{stheta}) gives
\begin{equation}\label{prodbound1}
 |\rho_m(\ell^TM\tau) - 1| \leq \prod_{j =1}^m \left(1 + \frac{|\ell_j|\pi}{2} |e^{2\pi i (M\tau)_j} -1| \right) - 1, \tau \in \mathbb R^n.
\end{equation}
$R^n$ acts on $\mathbb T^m$ to give a dynamical system
\begin{equation}\label{RnTm1}
 \pi^x(g) := \rho_m(Mx)g, \ \ x \in \mathbb R^n, \, g \in \mathbb T^m.
\end{equation}
This action is uniformly equicontinuous and since $\mathbb T^m$ is compact every orbit is almost periodic. Theregore for every $\delta > 0$ there exists a syndetic $S_M(\delta) \subset \mathbb R^n$ such that
\begin{equation}\label{RnTm2}
 \max_{j = 1,...,m} |e^{2\pi i (M\tau)_j} -1| \leq \delta, \ \ \tau \in S_M(\delta).
\end{equation}
For every $\epsilon > 0$ choose
\begin{equation}\label{epsdelta}
 \delta := \min \left\{ \frac{1}{m \pi} \left( \max_{\ell \in \Omega(q), j = 1,...,m} |\ell_j|\right)^{-1}, \,
\frac{\epsilon}{2mC_q} \right\}.
\end{equation}
(\ref{pbound1}), (\ref{prodbound1}) and the fact that $0 \leq \sigma \leq m/2$ implies $(1+\sigma)-1 \leq 2m\sigma$ gives
\begin{equation}\label{pbound2}
|p(\cdot + \tau) - p| \leq \epsilon, \ \ \tau \in S_M(\delta).
\end{equation}
\subsection{Minkowski Geometry and Mixed Volumes}\label{subsec2.8}
For any compact convex $K \subset \mathbb R^n,$ $V_n(K) := n$--dimensional volume of $K.$ The Minkowski sum of compact convex $K_i$ and $K_j$ is $K_i + K_j := \{x+y: x \in K_i, y \in K_j\},$
In \cite{bernshtein} Bernshtein defined the mixed volume of an $n$-tuple of compact convex sets
$V(K_1,...,K_n) :=$
$$
(-1)^{n-1} \sum_{i} V_n(K_i) + (-1)^{n-2} \sum_{i < j} V_n(K_i+K_j) + ... + V_n(K_1+\cdots+K_n).
$$
and noted that it differs from the classical definition by the factor $n!$ We will use Bernshtein's definition in the remainder of this paper.
\begin{lem} 
$V(K_1,...,K_n) \geq 0.$ 
$V(K_1,...,K_n) > 0$ iff there exists points $a_i, b_i \in K_i, i = 1,...,n$ such that $\{b_i-a_i: i = 1,...,n\}$ is linearly independent. Then
$V(K_1,...,K_n) \geq V([a_1,b_1],...,[a_n,b_n]) =  |\det [b_1-a_1, ...,b_n-a_n] |.$  
\end{lem}
Proof. (\cite{schneider}, Theorem 1.5.7, Theorem 5.1.8).
\begin{lem} 
If the vertices of $K_j$ are in $\mathbb Z^n$ then $V(K_1,...,K_n) \in \mathbb Z.$ 
\end{lem}
Proof. (\cite{ewald}, Theorem 3.9). 
\subsection{Genericity for Laurent Maps}\label{subsec2.9}
The fundamental theorem of algebra ensures that if $q \in \mathcal L_n$ then
$q$ has $V_1(\mathcal N(q))$ roots in $\mathcal C^*$ counting multiplicity. The situation for Laurent maps is more complicated even for a system of two linear equations in two variables. Consider $Q \in \mathcal L_2^2$ where $Q([z_1,z_2]) := [az_1 +  b z_2 - e,  cz_1+dz_2 - f]^T \in \mathcal L_2^2,$
$a, b, c, d, e, f \in \mathbb C^*.$ Then $\mathcal N(Q) = [K,K]^T$ where $K$ is the triangle with vertices
$[0\, 0]^T, [1\, 0]^T, [0\, 1]^T.$ A computation shows $|\Lambda(Q)| = 1$ iff 
$ad-bc \neq 0, de-bf \neq 0, -ce + af \neq 0$ and then $Q$ is called generic.
Otherwise either $\Lambda(Q) = \emptyset$ or  $\Lambda(Q)$ is unbounded and has no isolated points. We note that each of the three inequalities that characterize genericity involves amplitudes associated with vertices on one side of $K.$
\\ \\ 
For $Q \in \mathcal L_n^n,$ $\Lambda(Q)$ is analytic so is finite iff it is compact 
(\cite{chirka}, Proposition 1, p. 31) and genericity of $Q$ was 
defined by Bernshtein \cite{bernshtein} who proved 
\begin{lem}\label{bernshtein} 
If $Q$ is generic then $\Lambda(Q)$ is finite and
\begin{equation}\label{bern1}
	\sum_{\sigma \in \Lambda(Q)} m_\sigma(Q) = V(\mathcal N(Q)).
\end{equation}
If $Q$ is not generic then either $\Lambda(Q) = \emptyset$  or the set $\Lambda(Q)$ is noncompact and its subset of isolated roots 
$\Lambda_i(Q)$  satisfies
\begin{equation}\label{bern2}
	\sum_{\sigma \in \Lambda_i(Q)} m_\sigma(Q) 
< V(\mathcal N(Q)) \hbox{ if } V(\mathcal N(Q)) > 0,
\end{equation} 
and 
\begin{equation}\label{bern3}
	\sum_{\sigma \in \Lambda_i(Q)} m_\sigma(Q) 
= 0 \hbox{ if } V(\mathcal N(Q)) = 0.
\end{equation} 
\end{lem}
Fix $n \geq 1$ and for  $j = 1,...,n$ let $K_j$ be a polytope whose set of vertices $\mathcal V_j \subset \mathbb Z^n.$ 
The set of $Q \in \mathcal L_n^n$ with 
$\mathcal N(Q) = [K_1,...,K_n]^T$
is parameterized by
$$A := \mathbb C^{*k_1} \times \mathbb C^{k_2-k_1}$$ 
where $k_1 = \sum_{j=1}^n |\mathcal V_j|$ and $k_2 = \sum_{j=1}^n |\mathcal K_j \cap \mathbb Z^n|.$
\\ \\
We define genericity for Laurent maps. \\
For $u \in \mathbb R^n\backslash \{0\}$ and
$q = \sum_{\ell \in \Omega(q)} c_\ell \, \widetilde \ell \in \mathcal L_n$ define 
\\
$m(u,\Omega(q)) := \min \{u \cdot \ell : \ell \in \Omega(q)\},$ \\
$\Omega(q)_u := \{\ell \in \Omega(q) : u \cdot \ell = m(u,\Omega(q))\},$ \\
$q_u := \sum_{\ell \in \Omega(q)_u} c_\ell \, \widetilde \ell,$ 
and for  $Q = [q_1,...,q_n]^T \in \mathcal L_n^n,$ \\
$Q_u := [q_{1u},...,q_{nu}]^T.$
\begin{defi}\label{genericQ}
$Q \in \mathcal L_n^n$ is generic if $\Lambda(Q_u) = \emptyset$
for every $u \in \mathbb R^n \backslash \{0\}.$
\end{defi}
$\Omega(q)_u$ is a face of $\mathcal N(q)$ that is orthogonal to $u$ and has dimension $\leq n-1.$ Therefore there are only a finite number of $Q_u.$ Furthermore, Bernshtein (\cite{bernshtein}, p. 184) noted that each
$q_u$ "depends essentially on smaller number of variables" and therefore
"does not have any zeros in the general case." He showed that the set of non generic $Q$ are parameterized by a proper algebraic subset of the set of parameters $A$ described above.
Therefore if $Q$ is generic and $Q_1 \in \mathcal Ln^n$ satisfies 
$\mathcal N(Q_1) = \mathcal N(Q)$ and the amplitudes (coefficients) of $Q_1$ are sufficiently close to the amplitudes of $Q,$ then $Q_1$ is generic.
\\ \\
We illustrate genericity with an example more complex than the previous example. Consider the set of 
$Q = [q_1\ q_2]^T \in \mathcal L_2^2$ where 
$\mathcal N(Q) = [K \, K]^T$ and
$K$ is the square with vertices $[0 \, 0]^T, [1 \, 0]^T,[0 \, 1]^T,[1 \, 1]^T.$ These set of Laurent polynomials are parametrized by $C^{*8}.$
Then
$q_1([z_1,z_2]^T) := a_{0,0} + a_{1,0}z_1 + a_{0,1}z_2 + a_{1,1}z_1z_2,$
$q_2([z_1,z_2]^T) := b_{0,0} + b_{1,0}z_1 + b_{0,1}z_2 + b_{1,1}z_1z_2,$
where all $8$ coefficients are in $\mathbb C^*.$
If $u = r[\cos \theta, \sin \theta]^T$ and $r > 0,$ then 
$$
Q_u([z_1,z_2]^T) = \begin{cases}
[a_{1,0}z_1 + a_{1,1}z_1z_2, b_{1,0}z_1 + b_{1,1}z_1z_2]^T, \ \ \theta = 0 \\
[a_{0,1}z_2 + a_{1,1}z_1z_2, b_{0,1}z_2 + b_{1,1}z_1z_2]^T, \ \ \theta = \pi/2 \\
[a_{0,0} + a_{0,1}z_2, b_{0,0} + b_{0,1}z_2]^T, \ \ \theta = \pi \\
[a_{0,0} + a_{1,0}z_1, b_{0,0} + b_{1,0}z_1]^T, \ \ \theta = 3\pi/2 \\
[a_{1,1}z_1z_2,b_{1,1}z_1z_2]^T, \ \ \theta \in (0,\pi/2) \\
[a_{0,1}z_2,b_{0,1}z_2]^T, \ \ \theta \in (\pi/2,\pi) \\
[a_{0,0},b_{0,0}]^T, \ \ \theta \in (\pi,3\pi/2) \\
[a_{1,0}z_1,b_{1,0}z_1]^T, \ \ \theta \in (3\pi/2, \pi) 
\end{cases}
$$
$Q$ is generic when all $\Lambda(Q_u) = \emptyset.$
Clearly the last four Laurent maps above have no roots so a  calculation shows that $Q$ fails to be generic on the algebraic 
subset of $\mathbb C^{*4}$ 
that is the union of four hypersurfaces defined by 
$$
a_{1,0}\, b_{1,1} - a_{1,1}\, b_{1,0} = 0, \ \ \ \ 
a_{0,1}\, b_{1,1} - a_{1,1}\, b_{0,1} = 0,
$$
$$
a_{0,0}\, b_{0,1} - a_{0,1}\, b_{0,0}= 0, \ \ \ \ 
a_{0,0}\, b_{1,0} - a_{1,0}\, b_{0,0} = 0.
$$
If $Q$ is generic then calculation shows that $Q$ has two roots and the mixed volume $V(\mathcal N(Q)) = 2.$
\subsection{Genericity and Uniform Genericity for Trigonometric Maps}\label{subsec2.10}
Every nonzero $p \in \mathcal T_1$ has the form 
$p = c_1 \widehat \omega_1 + \cdots + c_d \widehat \omega_d$ 
where $d \geq 1,$ $c_j \in \mathbb C^*$ and
$\omega_1 < \cdots < \omega_d$ are real. Then for sufficiently large
$\Im z,$ $c_1 \widehat \omega_1(z)$ dominates the other terms in $p(z),$ and for
sufficiently large $-\Im z,$ $c_d \widehat \omega_d(z)$ dominates the other terms in $p(z).$ In both cases  $\Im \Lambda(p)$ is bounded. Moreover, $\Im \Lambda(p_2)$ is uniformly bounded for $p_1$ sufficiently close to $p$ and $\mathcal N(p_1)$ sufficiently close to $\mathcal N(p).$ Since $p$ is a nonzero univariate holomorphic function, $\Lambda(p) = \emptyset$ iff $d = 1$ otherwise it is an infinite discrete set. Corollary \ref{corProp3} implies $p$ is Bohr almost periodic. Hence corollary \ref{Rouchcor3}
implies $div \, p$ is Bohr almost periodic so (\ref{density1}) implies that 
its density $\Delta(div \, p)$ exists. A contour integration argument gives  
$\Delta(div \, p) = V(\mathcal N(p)).$ The situation for trigonometric maps is much more complicated and we must impose stringent conditions on 
$P \in \mathcal T_n^n$ to ensure that $\Im(\Lambda(P))$ is bounded, $\Lambda(P)$ is discrete and its density exists, and that these properties are preserved under sufficiently small perturbations.
\\ \\
We define genericity for trigonometric maps. \\
For $u \in \mathbb R^n\backslash \{0\}$ and
$p = \sum_{\omega \in \Omega(p)} c_\omega \, \widehat \omega \in \mathcal T_n$ define
\\
$m(u,\Omega(p)) := \min \{u \cdot \omega : \omega \in \Omega(q)\},$ \\
$\Omega(p)_u = \{\omega \in \Omega(p) : u \cdot \omega = m(u,\Omega(p))\},$ \\
$p_u := \sum_{\omega \in \Omega(p)_u} c_\omega \, \widehat \omega,$ and for $P = [p_1,...,n]^T \in \mathcal T_n^n$ define \\ 
$P_u := [p_{1u},...,p_{nu}]^T.$
\begin{defi}\label{genericP} 
$P \in \mathcal T_n^n$ is generic if $\Lambda(P_u) = \emptyset$ for every $u \in \mathbb R^n \backslash \{0\},$ 
\end{defi}
Let
$m := \hbox{dim }_{\mathbb Z} \Omega(P),$ $M \in \mathbb \mathbb R^{m \times n},$  $P := Q \circ \rho_m \circ M,$ and for $g \in \mathbb T^m$
define
\begin{equation}\label{Pg}
	P_g := Q \circ g\rho_m \circ M.
\end{equation}
We observe that if $\overline {\rho_m(M\mathbb R^n)} = \mathbb T^m$ 
then $\{P_g : g \in \mathbb T^m\} = \overline{O(P)}.$
We note that $\mathcal N(P_g) = \mathcal N(P)$ and only the arguments but not the moduli of the amplitudes of $P_g$ differ from those of $P.$
\begin{defi}\label{unifgenericP}
$P \in \mathcal T_n^n$ is uniformly generic if $P_g$ is generic for all $g \in \mathbb T^m.$
\end{defi}
Gelfond \cite{gelfond2} announced that if $P$ is uniformly generic then
$\Lambda(P)$ is discrete, $\Im(\Lambda(P))$ is bounded, and there exists 
$\beta > 0$ with
\begin{equation}\label{dendivP1}
	\left| 
\frac{r^{-n}}{c_n} \sum_{\lambda \in \Lambda(P), \Re \lambda \in B_n(0,r)}
m_\lambda(P) - V(\mathcal N(P)) 
\right| < \frac{\beta}{r}, \ \ r > 1.
\end{equation}
(\ref{dendivP1}) implies that the density $\Delta(div\, P) = V(\mathcal N(P)).$
\begin{remark}\label{Gelfond}
Gelfond deals with holomorphic Bohr almost periodic maps on vertical tube domains and counts roots whose imaginary parts lie in a cube $[-r/2,r/2]^n.$ He calls his condition properly unfolded. We avoid this term here because we reserve it for a condition that for every $u \in \mathbb R^n,$ $\Omega(p_j)_u$ consists of a single vertex for some $j.$ This condition implies uniform genericity but is stronger.
Gelfond's result applies to trigonometric maps since corollary \ref{Rouchcor3} implies they are Bohr almost periodic. For trigonometric maps the indicator diagrams of the components are their Newton polytopes so Gelfond's condition (\cite{gelfond1}, System 3, p. 155) applies and the requirement that $P$ be uniformly generic is stated in (\cite{gelfond1}, (c), p. 156) "it is required that this closure does not contain consistent systems" meaning $u \in \mathbb R^n\backslash \{0\}$ and $P_u \neq \emptyset.$ Gelfond's paper does not provide a proof but references his preprint \cite{gelfond0}. He also references
a paper of Hovanskii \cite{hovanskii2}, in which 4 out of 7 references are preprints, and a paper by Kazarnovskii \cite{kazarnovskii1} in Russian. Kazarnovskii's paper \cite{kazarnovskii2} derived asymptotic distribution of roots for exponential sums whose exponents are not necessarily purely imaginary and thus not trigonometric polynomials. 
\end{remark}
\begin{prop}\label{prop4} 
If $P \in \mathcal T_n^n$ is uniformly generic then there exists a tube domain $T$ such that for all $P_1 \in \mathcal T_n^n$
sufficiently close to $P,$ $P_1 \in B(T)$ and $div\, P_1$ is a Bohr almost periodic multiset contained in $T.$ 
Moreover, $P_1$ is real-rooted iff $\Delta(\Lambda(P_1) \cap \mathbb R^n) = \Delta(\Lambda(P_1)).$
\end{prop}
Proof. Follows from Corollaries 1, 3, and 4. 
\\ \\
The following result extends Levin's result (\cite{levin}, p. 268, Lemma 1) for $n = 1.$
\begin{prop}\label{prop5} 
If $P \in \mathcal T_n^n$ is uniformly generic and $T$ is a tube domain containing $\Lambda(P),$
$\epsilon > 0,$ $T_\epsilon$ is the set of $z \in T$  
such that $|z-w| \geq \epsilon$ whenever $P(w) = 0$ or
$w \notin T.$ Then $\inf_{z \in T_\epsilon} |P(z)| > 0.$
\end{prop}
Proof. Follows from Corollaries 2 and 4.
\subsection{Pontryagin Duality for Locally Compact Abelian Groups}\label{subsec2.11} 
Let $H$ be a locally compact abelian group. A continuous homomorphism $\chi : H \mapsto \mathbb T$ is called a character of $H.$ The Pontryagin dual $\widehat H$ of $H$ is the set of characters of $H$ equipped with pointwise multiplication and the topology of uniform convergence on compact subsets. $\widehat H$ is a locally compact abelian group with identity the constant function with value $1.$ We observe that 
$$\widehat {\mathbb T^m} = \left\{\widetilde \ell|_{\mathbb T^m} : \ell \in \mathbb Z^m \right\} \cong \mathbb Z^m.$$ 
$$\widehat {\mathbb R^n} = \left\{ \widehat \omega|_{\mathbb R^n} : \omega \in \mathbb R^n \right\} \cong \mathbb R^n.$$
$H$ is discrete (resp. compact) iff $\widehat H$ is compact (resp. discrete).
If $H$ is the additive group of a field then $\widehat H$ is isomorphic to $H,$ though not in a unique way. $H$ is connected iff $\widehat H$ is torsion free (no nonidentity elements has finite order). $\hbox{dim }H = \hbox{rank }\widehat H.$ 
If $H_1$ and $H_2$ are locally compact topological groups and $\varphi : H_1 \mapsto H_2$ is a continuous homomorphism
then its Pontryagin dual $\widehat \varphi : \widehat H_2 \mapsto \widehat H_1,$ defined by 
$$
\widehat \varphi(\chi) := \chi \circ \varphi, \ \ \chi \in \widehat H_2,
$$ 
is a continuous homomorphism. Moreover, 
$\varphi$ is injective (resp. has dense image) iff $\widehat \varphi$ has dense image
(resp. is injective).
Pontryagin's duality theorem (\cite{rudin1}, Theorem 1.7.2) implies that the inclusion 
$\iota : H \mapsto \widehat {\widehat H},$ defined by
$$
\iota(h)(\chi) := \chi(h), \ \chi \in \widehat H,
$$ 
is a isomorphism of topological groups, so
we identify  $\widehat {\widehat H}$ with $H$ and if $\varphi : H_1 \mapsto H_2$ we identity
$\widehat {\widehat \varphi}$ with $\varphi.$
Therefore characters on $H$ separate points of $H$ and
$H \neq \{1\}$ iff $\widehat H \neq \{1\}.$
If $H \subset G$ is a closed subgroup of a locally compact abelian group $G$ define
$$
H^\perp := \{\, \chi \in \widehat G \, : \, \chi(H) = \{1\} \, \}.
$$
Pontryagin duality gives 
$H = \{ \, h \in G \, : \, \chi(h) = 1 \hbox{ for every } \chi \in H^\perp \, \}.$
If $G$ is a closed subgroup of $H$ then the quotient $G/H$ is a locally compact abelian group and the map
$\tau : H \mapsto H/G$ defined by $\tau(h) := hG$ is a continuous surjective homomorphism with kernel $G.$
Hence $G \neq H$ iff $H/G \neq \{1\}$ iff there exists $\chi \in \widehat {H/G} \backslash \{1\}$ iff $\chi \circ \tau \in \widehat H \backslash \{1\}$ with $\chi(G) = \{1\}.$ 
Therefore if $G$ is a closed subgroup of $\mathbb T^m,$ then $G \neq \mathbb T^m$ iff there exists 
$\ell \in \mathbb Z^m \backslash \{0\}$ such that $\widetilde \ell(G) = \{1\},$ and if
$G$ is a closed subgroup of $\mathbb R^n,$ then $G \neq \mathbb R^n$ iff there exists 
$w \in \mathbb R^n \backslash \{0\}$ such that $\widehat w(G) = \{1\}.$
\\ \\
The following closure results are used throughhout the sequel of this paper.
\begin{lem}\label{closure} 
Assume $M \in \mathbb R^{m \times n},$
$u \in \mathbb R^m,$ $\psi := \rho_m \circ M : \mathbb R^n \mapsto \mathbb T^m,$ and
$\varphi : \mathbb Z \mapsto \mathbb T^m$  is defined by $\varphi(k) := \rho_m(ku).$ Then
\begin{enumerate}
\item $\overline {\psi(\mathbb R^n)} = \mathbb T^m \text{ iff } \text{dim}_{\, \mathbb Z} M = m.$
\item $\text{kernel } \psi = \{0\} \text{ iff rank } M = n \text{ and } \mathbb Z^m \cap M\mathbb R^n = \{0\}.$ 
\item $\overline {M^T\mathbb Z^m} = \mathbb R^n$ iff $\text{kernel } \psi = \{0\}.$  
\item $\overline {\varphi(\mathbb Z)}= \mathbb T^m$ iff
$\text{dim}_{\, \mathbb Z} 
\{1,u_1,...,u_m\} = m+1.$
\end{enumerate}
\end{lem}
Proof. 1. follows from (\ref{hatw2}) since 
$\overline {\psi(R^n)} \neq \mathbb T^m$ iff there exists $\ell \in \mathbb Z^m \backslash \{0\}$ with 
$$\{1\} = \widetilde \ell\, (\psi(\mathbb R^n)) = 
\widehat {M^T\ell}\, (\mathbb R^n) \iff
 M^T\ell = 0$$
and the last inclusion holds for some $\ell \in \mathbb Z^m \backslash \{0\}$ iff
$\text{dim}_{\, \mathbb Z} M < m.$
\\
2. follows since kernel $\rho_m = \mathbb Z^m$ implies
$\text{kernel } \psi = 
\{x \in \mathbb R^n : Mx \in \mathbb Z^m\}$ so
$\{0\} = \text{ kernel } \psi $
iff $\text{ rank } M = n \text{ and } \mathbb Z^m \cap M\mathbb R^n = \{0\}.$
\\
3. follows since $\overline {M^T\mathbb Z^m} \neq \mathbb R^n$ iiff there exists $\omega \in \mathbb R^m\backslash \{0\}$ with
$$
\{1\} = \widehat \omega \, (M^T\mathbb Z^m) = \rho_1((M\omega)\cdot \mathbb Z^m) \iff M\omega \in \mathbb Z^m$$
and the last equality holds for some $\omega \in \mathbb R^n \backslash \{0\}$ iff
$\text{kernel } \psi \neq \{0\} .$
\\
4. follows since $\overline {\varphi(\mathbb Z)}= \mathbb T^m$
iff there exists $\ell \in \mathbb Z^m \backslash \{0\}$ with 
$$\{1\} = \widetilde \ell\, (\varphi(\mathbb Z)) = \widetilde \ell\, (\rho_m(\mathbb Zu)) = \rho_1(\ell \cdot u)$$
iff $\ell \cdot u \in \mathbb Z$ iff 
$\text{dim}_{\, \mathbb Z} 
\{1,u_1,...,u_m\} \le m.$
\subsection{Compactification of Groups, Functions, Measures. Multisets}\label{subsec2.12}
\begin{defi}\label{compactification}
A compactification of a locally compact abelian group $H$ is a pair 
$(G,\psi)$ where $G$ is a compact group and $\psi : H \mapsto G$ is a continuous homomorphism with a dense image. Clearly $G$ is abelian and if $H$ is connected then $G$ is connected. The Bohr compactification $(G_1,\psi_1)$ is when
$G_1 = \widehat {\widehat H_d}$ where $\widehat H_d$ equals $\widehat H$ with the discrete topology and $\psi_1 = \widehat \iota$ and
$\iota :  \widehat H_d \mapsto \widehat H$ is the identity homomorphism.
\end{defi}
The Bohr compactification has the following universal property: If $(G,\psi)$ is any compactification of $H$ then there exists a surjective homomorphism
$\Psi : G_1 \mapsto G$ such that $\psi = \Psi \circ \psi_1.$
If $(G,\psi)$ is a compactification of $H$ we define its spectrum
$\Omega(G,\psi) := \widehat \psi (\widehat G).$ 
If $D$ is the group generated by $\Omega(G,\psi)$ with the discrete topology then
Pontryagin duality gives $G =\widehat D.$
Let $\mathcal T(G) :=$ the algebra of finite linear 
combinations of characters on $G.$
Since $\mathcal T(G)$ separates points and is closed under complex conjugation, 
the Stone-Weierstrass theorem (\cite{rudin1}, A14) implies that it is dense in $C(G).$
For $f \in B(\mathbb R^n)$ let $D_f :=$ the discrete subgroup of $\mathbb R^n$ 
generated by $\Omega(f).$
\begin{lem}\label{fF} 
If $F \in C(G)$ and $f := F \circ \psi$ then $f \in B(\mathbb R^n)$ and $\mathcal M(f) = \int_G F$
where integration is with respect to Haar measure on $G.$
Let  $\varphi : D_f \mapsto \mathbb R^n$ be the inclusion map.
If $f \in B(\mathbb R^n)$ then $(G,\psi),$ where $G := \widehat D_f$ and $\psi := \widehat \varphi,$ 
is a compactification of $\mathbb R^n.$ Furthermore $\Omega(f) \subset \Omega(G,\psi)$ and there 
exists a unique $F \in C(G)$ such that $f = F \circ \psi.$
\end{lem}
Proof. Since $\mathcal T(G)$ is dense in $C(G)$ $F$ can be uniformly approximated $\tau \mathcal T(G)$ hence
$p$ can be approximated by $\tau \circ \psi \in \mathcal T_n$ so $f \in B(\mathbb R^n).$ The integral identity
occurs because $\int_G \chi = 1$ if $\chi = 1$ else $\int_G \chi = 0.$ 
Since $D_f$ is discrete $G$ is compact and since $\varphi$ is injective $\psi$ has a dense image hence
$(G,\psi)$ is a compactification of $\mathbb R^n.$ Pontryagin duality implies that
$\widehat G = D_f$ and $\widehat \psi = \varphi$ hence $\Omega(G,\psi) = D_f.$
Then $\Omega(f) \subset \Omega(G,\psi)$ so there exist a sequence
$p_k \in \mathcal T_n$ with $\Omega(p_k) \subset \Omega(G,\psi)$ and $||f-p_k||_\infty \rightarrow 0.$
Since $\Omega(p_k) \subset \Omega(G,\psi)$ there exists a sequence $\tau_k \in \mathcal T(G)$ such that
$p_k = \tau_k \circ \psi.$  Since $\psi$ has dense image $||\tau_j-\tau_k||_\infty = ||p_j-p_k||_\infty.$
Since $p_k$ converges it is a Cauchy sequence hence $\tau_k$ is a Cauchy sequence so converges to a function
$F \in C(G).$ Then $f = F \circ \psi.$
\begin{defi}
Let $\mathcal B(G)$ denote the set of Borel measures on a compact abelian group $G.$ 
The Riesz Representation theorem (\cite{rudin1}, Theorem 6.19) implies that 
$\mathcal B(G)$ is the space of continuous linear functionals on $C(G).$ For $\widetilde \mu \in \mathcal B(G)$
we define its Fourier transform $\widehat {\widetilde \mu} : \widehat G \mapsto \mathbb C$
by 
$$\widehat {\widetilde \mu}(\chi) := \, <\widetilde \mu , \, \overline \chi>\, , \  \ \chi \in \widehat G,$$
and its spectrum
$$\Omega(\widetilde \mu) := \{\chi \in \widehat G : \widehat {\widetilde \mu}(\chi) \neq 0 \}.$$
\end{defi}
\begin{defi}\label{approxid}
An approximation of the identity on $\mathbb R^n$ is a sequence
$\xi_j \in C_c(\mathbb R^n)$ that is nonnegative, $\int_{x \in \mathbb R^n}\xi_j(x)dx = 1,$ and for every open set $U$ containing $0,$ there exists $J(U) \in Z_+$ such that $j \geq N(U) \implies$ support $\xi_j \subset U.$ Then $\xi_k$ converges weakly to $\delta_0.$
\end{defi}
%
%
\begin{prop}\label{prop6}
If $\mu$ is a Bohr almost periodic measure on $\mathbb R^n,$
then there exists a compactification $(G,\psi)$ of $\mathbb R^n$ and
$\widetilde \mu \in \mathcal B(G)$ such that
\begin{equation}\label{hatmu}
\mathcal F_B(\mu)(\chi\circ\psi) = \widehat {\widetilde \mu}(\chi), \ \ \chi \in \widehat G. 
\end{equation}
The measure $\widetilde \mu$ is called the trace of $\mu$
with respect to $(G,\psi)$ and it is unique.
\end{prop}
Proof. Let $D_\mu$ be the discrete subgroup of $\mathbb R^n$ generated by $\Omega(\mu)$ 
then construct the compactification $(G,\psi)$ of $\mathbb R^n$ as in the proof of Lemma \ref{fF}. 
Let $\xi_j$ be an approximation of the identity on $\mathbb R^n.$ Since $\mu$ is translation bounded
\begin{equation}\label{cond}
 \text{sup}_j \, \mathcal M(\xi_j*\mu) < \infty,
\end{equation}
Since $\mu$ is a Bohr almost periodic measure 
there exists a sequence $\Xi_j \in C(G)$ of such that
$\xi_j = \Xi \circ \psi.$ Furthermore $\Xi_j$ are nonnegative and
$\mathcal M(\xi_j*\mu) = \int_G \Xi_j$ so
$\text{sup} \int_G \Xi_j < \infty.$ 
Therefore the Banach-Alaoglu Theorem (\cite{rudin3}, 3.15) implies that
there exists a subsequence of $\Xi_j$ that converges weakly to a measure $\widetilde \mu$ on $G.$ 
Replacing the sequences with this subsequence and applying 
(\ref{FBconvf}) and Lemma \ref{fF} gives
\begin{equation}
\begin{array}{ccc}
\mathcal F_B(\mu)(\chi \circ \psi) & = & \lim_{j \rightarrow \infty} \widehat \xi_j (\chi \circ \psi) \mathcal F_B(\mu)(\chi \circ \psi) \\
= \lim_{j \rightarrow \infty} \mathcal F_B(\xi_j*\mu)(\chi \circ \psi) & = & \lim_{j \rightarrow \infty} \mathcal F_B(\Xi_j \circ \psi)(\chi \circ \psi) \\
= \lim_{j \rightarrow \infty} \mathcal M((\Xi_j \circ \psi)(\overline \chi \circ \psi)) & = & \lim_{j \rightarrow \infty} \mathcal M((\Xi_j \,\overline \chi) \circ \psi) \\
= \lim_{j \rightarrow \infty} \int_G \Xi_j \overline \chi & = & \int_G \overline \chi \, d \widetilde \mu. 
\end{array}
\end{equation}
Since $\mathcal T(G)$ is dense in $C(G),$ 
$\widetilde \mu$ is uniquely determined $\widehat {\widetilde \mu}.$ 
\begin{remark}
The trace of $\mu$ with respect to the Bohr compactification equals the image of $\mu$ under the Bohr map  (\cite{argabright}, Section 7).
\end{remark}
\section{Examples}\label{sec3} 
Example 1 constructs a trivial Fourier quasicrystal that is not uniformly discrete, hence not Delone. Example 2 constructs two model sets that are Besicovitch but not Bohr almost periodic. Moreover, neither of their associated measures are Poisson because their Fourier-Bohr transforms do not satisfy (\ref{summable}). However, their sum is a Bohr almost periodic Poisson measure. Example 3 constructs a family of uniformly discrete Bohr almost periodic subsets whose associated measures are Fourier quasicrystals. The algebraic varieties they use are generally not complete intersections (codimension $n$ and defined by the common zeros of $n$ Laurent polynomials). In this respect they differ from the construction of Fourier quasicrystals in Section 5 but are similar to the construction by Alon, Kummer, Kurasov and  Vinzant \cite{alon2}.
\subsection{Almost Periodic Multisets}\label{subsec4.1}
\begin{example}\label{aZbZ} 
Choose $a, b \in \mathbb R_+,$ with $a/b \notin \mathbb Q,$ 
and define three multisets
$$\alpha := (a\mathbb Z,m_\alpha), m_\alpha := 1 \text{ on }a\mathbb Z,
\beta := (b\mathbb Z,m_\beta), m_\beta := 1 \text{ on } b\mathbb Z,
\gamma := \alpha \oplus \beta
$$ 
where $\oplus$ is explained in Subsection \ref{subsec2.5}.
\end{example}
We record the following observations without proof: these are FQs on $\mathbb R$ that by Definition \ref{trivial} are trivial, their indices defined by (\ref{ind}) satisfy $ind(\alpha) = ind(\beta) = 1, ind(\gamma) = 2.$ Since these multisets are Bohr almost periodic, their orbit closures are compact groups. Furthermore, computations give  
$$
\begin{array}{cc}
\overline {O(\alpha)} = & \{\pi^x(\alpha) : x \in \mathbb R\} \cong \mathbb T, \\
\overline {O(\beta)} = & \{\pi^x(\beta)  : x \in \mathbb R\} \cong \mathbb T, \\
\overline {O(\gamma)} = & \{\pi^x(\alpha) \oplus \pi^y(\beta) : x, y \in \mathbb R^n\} \cong \mathbb T^2.
\end{array}
$$
\begin{example}\label{alphaj} 
Parametrize $\mathbb T^2$ by $[0,1)^2$ with coordinates $x, y \in [0,1).$ Let $\theta = \pi/3$ and 
define the homomorphism $\psi : \mathbb R \mapsto \mathbb T^2$ by
$$
	\psi(t) := t[\cos \theta, \sin \theta]^T \mod \mathbb Z^2.
$$
Since $\tan \theta = \sqrt 3$ is irrational, $\psi$ has a dense image so $(\mathbb T^2,\psi)$ is a compactification of $\mathbb R.$
Let $c := \sqrt 2$  and define three $1$-chains in $\mathbb T^2$
$$
\begin{array}{cc}
L_1 := & \left \{ 
\left[\begin{array}{c}
x \\
y
\end{array}\right]
: x \in [0,1/2), y = cx\right \}, \\
L_2 := & \left \{ 
\left[\begin{array}{c}
x \\
y
\end{array}\right]
: x \in [1/2,1), y = c-cx\right \}, \\
L_3 := & L_1 \cup L_2.
\end{array}
$$
Since $0 < c = \sqrt 2 < \tan \theta = \sqrt 3,$ each $L_j$ is transverse to the foliation on $\mathbb T^2$ whose leaves are the cosets of 
$\psi(\mathbb R).$ Therefore each $\psi^{-1}(L_j)$ is a uniformly discrete subset of $\mathbb R$ so we can define the multiset $\alpha_j$ 
and associated pure point Radon measure $\mu_j$ by
$$
\alpha_j := (\psi^{-1}(L_j),m_j), \ m_j = 1 \text{ on } \psi^{-1}(L_j)
$$
$$
\mu_j := \sum_{\lambda \in \psi^{-1}(L_j)} \delta_\lambda.
$$
Clearly $\alpha_3 = \alpha_1 \oplus \alpha_2$ and 
$\mu_3 = \mu_1 + \mu_2.$
\end{example}
$\alpha_1$ and $\alpha_2$ are Meyer model sets \cite{moody}, (\cite{lawton}, Definition 2). Therefore, as noted in Subsection 
\ref{subsec2.4} they are Besicovitch almost periodic but not Bohr 
almost periodic. However, since $L_3$ is a cycle (closed loop) 
$\alpha_3$ is Bohr almost periodic. We will prove that
for every $a, b \in \mathbb R$ with $0 < b - a < 2/3,$
\begin{equation}\label{sumFmuj}
 \sum_{\omega \in (a,b)}
|\mathcal F_B(\mu_j)(\omega)| \ \ \ 
\begin{cases}
 \ \ = \infty, \ \ j = 1, 2, \\
\ \ < \infty, \ \ j = 3.
\end{cases}.
\end{equation}
Lemma \ref{meyer} and (\ref{sumFmuj}) imply that $\widehat \mu_1$ and $\widehat \mu_2$ are not Radon measures and $\widehat \mu_3$ is a Radon measure. Therefore $\mu_1$ and $\mu_2$ are not Poisson measures and $\mu_3$ is a Poisson measure. 
Our proof  combines results in \cite{lawton} about Bohr almost periodic sets and results in \cite{KN} about uniform distribution.  
\\ \\
(\cite{lawton}, Remark 4) implies that
$\Omega(\mu_j) \subset \mathbb Z \cos \theta + \mathbb Z \sin \theta, \ j = 1, 2, 3$
and for $\omega = \ell_1 \cos \theta + \ell_2 \sin \theta),$
\begin{equation}\label{FBmuj}
\mathcal F_B(\mu_j)(\omega) = 
\int_{L_j} e^{2\pi i (\ell_1x + \ell_2y)}\, (\sin \theta \ dx - \cos \theta \ dy), \ \  j = 1, 2, 3
\end{equation}
hence
\begin{equation}\label{eq69}
	\mathcal F_B(\mu_1)(\omega) = 
\frac{\sqrt 3 - \sqrt 2}{2\pi(\ell_1+\ell_2\sqrt 2 )} e^{-\pi i (\ell_1+\ell_2\sqrt 2 )/2}\, \sin \pi(\ell_1+\ell_2\sqrt 2)/2.
\end{equation}
\begin{equation}\label{eq70}
	\mathcal F_B(\mu_2)(\omega) = 
\frac{\sqrt 3 + \sqrt 2}{2\pi(\ell_1-\ell_2\sqrt 2 )}
\,e^{-\pi i (3\ell_1+\ell_2\sqrt 2)/2}\, \sin \pi(\ell_1-\ell_2\sqrt 2)/2.
\end{equation}
Since the quantities in (\ref{eq69}) and (\ref{eq70}) never vanish, \\
$\Omega(\mu_j) = 
(\mathbb Z \cos \theta + \mathbb Z \sin \theta)\backslash \{0\}, \ j = 1, 2$ 
which is dense in $\mathbb R$ by Lemma \ref{closure}.
For $\omega \in (a,b),$ 
$\ell_1 = -\ell_2 \sqrt 3 + 2g(\ell_2)$ where $g(\ell_2) \in (a,b),$ hence
\begin{equation}\label{mFBmu1}
|2\pi \ell_2 \mathcal F_B(\mu_1)(\omega)| = O(1/\ell_2) + \left|\sin \pi \left[\ell_2\frac{\sqrt 3 - \sqrt 2}{2} - g(\ell_2)\right]\right|.
\end{equation}
\begin{equation}\label{mFBmu2}
|2\pi \ell_2 \mathcal F_B(\mu_2)(\omega)| = O(1/\ell_2) +   \left|\sin \pi\left[\ell_2\frac{\sqrt 3 + \sqrt 2}{2} - g(\ell_2)\right]\right|.
\end{equation}
\begin{equation}\label{mFBmu3}
|2\pi \ell_2 \mathcal F_B(\mu_3)(\omega)| = O(1/\ell_2) \, + 
\end{equation} 
$$
\begin{cases}
\left|\sin \pi \left[\ell_2\frac{\sqrt 3}{2} - g(\ell_2)\right]
\cos \pi \left[\ell_2\frac{\sqrt 2}{2}\right]\right|, \ \ \ell_1 \text{ even}, \\ \\
\left|\cos \pi \left[\ell_2\frac{\sqrt 3}{2} - g(\ell_2)\right]
\sin \pi \left[\ell_2\frac{\sqrt 2}{2}\right]\right|, \ \ \ell_1 \text{ odd}
\end{cases} 
$$
Since $\ell_2\frac{\sqrt 3}{2} - g(\ell_2) = -\frac{\ell_1}{2},$
the second term in the right side of (\ref{mFBmu3}) equals $0,$ hence (\ref{sumFmuj}) holds for $\mu_3.$
A calculation shows that the first term in the right side of (\ref{mFBmu3})  never vanishes hence
$\Omega(\mu_3)$ is dense in $\mathbb R.$
\\ \\
Define $\varphi_1 : \mathbb Z \rightarrow \mathbb T^2$ by
$$
	\varphi_1(\ell_2) := \ell_2[\sqrt 3, (\sqrt 3 - \sqrt 2)/2\} = \text { mod } \mathbb Z^2.
$$
Since $dim_{\mathbb Z}\{1, \sqrt 3, (\sqrt 3 - \sqrt 2)/2\} =  3,$
Lemma \ref{closure} implies 
$\overline {\varphi_1(\mathbb Z^2}) = \mathbb T^2$ and
(\cite{KN}, Example 6.7) shows that the sequence 
$\varphi_1(\ell_2), \ell_2 \in \mathbb Z_+$ is uniformly distributed. This means that 
for every nonempty open $O \subset \mathbb T^2,$  
\begin{equation}\label{density1}
	\lim_{L \rightarrow \infty} L^{-1} 
|\\, [1,L] \cap\varphi_1^{-1}(O) \, | = h(O)
\end{equation}
where $h(O) := $Haar measure of $O.$ Since $O$ is non empty and open, $h(O) > 0$ hence the density of the sequence is positive. We observe that (\ref{density1}) implies 
\begin{equation}\label{infsum}
	\lim_{L \rightarrow \infty} 
\sum_{\ell_2 \in \varphi^{-1}(O) \cap [1,L]} \frac{1}{\ell_2} = \infty.
\end{equation}
Choosing 
$$O_1 :=  (2a,2b) \times
	(b+1/6,a+5/6)
\text{ mod } \mathbb Z^2$$
gives
$\ell_2 \in \varphi^{-1}(O_1)$ iff 
there exists $\ell_1 \in \mathbb Z$ with
$$\omega := \ell_1 \cos \theta + \ell_2 \sin \theta \in (a,b)$$ 
and
$$\ell_2 \frac{\sqrt 3 - \sqrt 2}{2} - g(\ell_2) \text{ mod } \mathbb Z  \in (1/6, 5/6).$$ Therefore the second term on the right side of
(\ref{mFBmu1}) is $> 1/2,$ hence (\ref{sumFmuj}) holds for $j = 1.$ The proof for $j = 2$ is similar. 
\subsection{Fourier Quasicrystals}\label{FQ}
We use results in \cite{lawton} to construct Fourier quasicrystals on 
$\mathbb R^n, n \geq 1.$ Let $\Lambda \subset \mathbb R^n$ be 
a uniformly discrete Bohr almost periodic set of toral type and 
$\mu := \sum_{\lambda \in \Lambda} \delta_\lambda$ be its 
associated Bohr almost periodic measure. 
(The abstract of \cite{lawton} says incorrectly that $\Lambda$ has toral type 
if the Fourier transform of $\mu$ equals zero outside of: a rank 
$m < \infty$ subgroup of $\mathbb R^n.$ It should say: a generated 
by $m < \infty$ elements subgroup of $\mathbb R^n.$ The 
group of dyadic rational numbers has rank $1$ but is not finitely 
generated.) The results in \cite{lawton} imply: 
\begin{enumerate}
\item the associated compactification of $\mathbb R^n$ 
equals $(\mathbb T^m,\psi),$ where $m \geq n.$ 
$m = n$ iff $\Lambda$ is a lattice subgroup of 
$\mathbb R^n.$ $\psi = \rho_m \circ M$ 
where $M \in \mathbb R^{m \times n}.$ 
\item $\Lambda = \psi^{-1}(X)$ where $X := \overline {\psi(\Lambda)}$ 
is an $(m-n)$-dimensional submanifold of 
$\mathbb T^m$ transverse to the foliation of 
$\mathbb T^m$ whose leaves are the cosets of $\psi(\mathbb R^n).$
Each connected component of $X$ is homotopic to a 
subgroup of $\mathbb T^m$ isomorpic to $T^{m-n},$ 
and the homotopy is transverse to the  foliation.
\item the Fourier transform $\widehat \mu$ is zero outside of $M^T\mathbb Z^m,$ and
\begin{equation}\label{FTmu}
\widehat \mu(M\ell) = \int_{z \in X} z^{-\ell} \, d\theta, \ \ \ell \in \mathbb Z^m
\end{equation}
where $d\theta$ is the unique differential $(m-n)$-form transverse to the folliation and
$\int_X d\theta = \Delta(\Lambda).$  The density $\Delta(\Lambda)$ is determined by $M$ and the homotopy classes of the connected components of $\overline {\psi(\Lambda)}.$
\end{enumerate}
These results imply that if $X$ has $J$ connected components, then $\Lambda$ is a union of $J$ irreducible Bohr almost periodic sets (meaning that they are not the union of Bohr alost periodic subsets). Furthermore, since $|\widehat \mu| \leq \Delta(\Lambda),$ (\ref{FTmu}) implies that $\mu,$ and hence $\Lambda,$ is a Fourier quasicrystall iff $|\Omega(\mu) \cap [-r,r]^n|$ is bounded by a polynomial in $r > 0.$
\\ \\
In this manuscript we use a different notation from the notation used in \cite{lawton}. 
We denote the Fourier transform of a tempered distribution
$\mu$ by $\widehat \mu,$ and denote the Fourier-Bohr transform of an almost periodic measure $\mu$ by $\mathcal F_B(\mu).$
\cite{lawton} denoted the Fourier transform of a tempered distribution
$\mu$ by $\mathcal F(\mu),$ and denoted the Fourier-Bohr transform of an almost periodic measure $\mu$ by $\widehat \mu.$
We define the Fourier-Bohr transform of an almost periodic set (or multiset) to be the Fourier-Bohr transform of its associated measure. Therefore, if $\mu_t$ is the measure corresponding to $\Lambda_t$ in the example below, then
$\mathcal F_B(\Lambda_t)(M^T\ell)$ in (\ref{FBmu1}) equals 
$\widehat {\mu_t}(M\ell)$ in the notation used in \cite{lawton}.
\begin{example}\label{example1}
Fix $n \in \mathbb Z_+,$ 
$b \in \mathbb R_+^n$ with
$dim_{\, \mathbb Z}\{1, b_1, ...,b_n\} = n+1,$
distinct $s_1,...,s_{n+1} \in (-1,1),$
$\gamma \in Z_+^{n+1}$ with $\hbox{gcd }(\gamma_1,...,\gamma_{n+1}) = 1,$ and define
\begin{equation}\label{mMpsi}
m :- n+1, \ \ 
M := \left[
\begin{array}{c}
I_n \\
-b^T 
\end{array}
\right] \in \mathbb R^{m \times n}, \ \ 
\psi := \rho_m \circ M : \mathbb R^n \mapsto \mathbb T^m.
\end{equation}
For every $t \in [0,1],$ 
define $\varphi_t : \mathbb T \rightarrow \mathbb T^m$ by
\begin{equation}\label{varphit}
\varphi_t(\zeta)_j := \frac{\zeta^{\gamma_j}-ts_j}{1-ts_j\zeta^{\gamma_j}}, \ \ \
\zeta \in \mathbb T, \  j = 1,...,m,
\end{equation}
$X_t := \varphi_t(\mathbb T),$ and
\begin{equation}\label{Lambdatdef}
\Lambda_t := \psi^{-1}(X_t).
\end{equation}
\end{example}
The assumptions on $b,$ definitions of $M$ and $\psi,$
and Lemma \ref{closure} imply that $(\mathbb T^m,\psi)$ is a compactification of $\mathbb R^n,$ $kernel(\psi) = \{0\},$ and $\overline {M\mathbb Z^m} = \mathbb R^n.$
\begin{lem}\label{transverse} 
Every $\varphi_t : \mathbb T \mapsto \mathbb T^m$ is injective.
$X_0$ is a subgroup of $\mathbb T^m$ transverse to the foliation on $\mathbb T^m.$ $\Lambda_0$ is a lattice subgroup of $\mathbb R^n$  whose density $\Delta(\Lambda_0)$ is determined by $M$ and $\gamma.$
The family of maps 
$\varphi_t : \mathbb T \mapsto \mathbb T^m,  t \in [0,1]$ is a homotopy so these maps induce identical maps on fundamental groups:
$$\pi_1(\varphi_t) : \pi_1(\mathbb T) \mapsto \pi_1(\mathbb T^m),$$
namely the map $\mathbb Z \mapsto \mathbb Z^m$ defined by
$\ell \mapsto \ell \gamma.$
There exists $\tau \in (0,1]$ such
that for every $t \in [0,\tau],$ $X_t$ 
is transverse to the foliation on $\mathbb T^m$ (whose leaves are the cosets of $\psi(\mathbb R^n),$
$\Lambda_t$ is a 
uniformly discrete Bohr almost periodic set that satisfies
$\overline {\psi(\Lambda_t)} = X_t,$ and is a Bohr almost periodic perturbation of $\Lambda_0,$ hence its density
$\Delta(\Lambda_t) = \Delta(\Lambda_0).$ The spectrum 
$\Omega(\Lambda_t) \subset M^T\mathbb Z^m.$ 
$\Lambda_t$ is irreducible. $\Lambda_t$ is nontrivial as defined by
Definition \ref{trivial}. 
\end{lem}
Proof. $gcd(\gamma_1,...,\gamma_n) = 1$ implies there exists $\ell \in \mathbb Z^m$ with
$\gamma \cdot \ell = 1.$ Furthermore (\ref{varphit}) implies
\begin{equation}\label{invvarphit}
\zeta^{\gamma_j} = \frac{\varphi_t(\zeta)_j+ts_j}{1+ts_j\varphi_t(\zeta)_j}, \ \ \zeta \in \mathbb T, \, j = 1,...,m.
\end{equation}
Therefore 
\begin{equation}\label{zeta}
\zeta = \zeta^{\gamma \cdot \ell} = 
\sum_{j=1}^m (\varphi_t(\zeta))^{\ell_j},
\end{equation}
hence $\varphi_t(\zeta)$ determines $\zeta$ so it
is injective. Therefore $X_0$ is a subgroup of $\mathbb T^m,$
$\varphi_0$ is an isomorphism of $\mathbb T$ onto $X_0,$ and $X_o$ is transverse to the foliation on $\mathbb T^m$ because 
$\gamma \notin M\mathbb R^n.$ 
Since $\varphi_t$ is continuously differentiable 
function of $t \in [0,1],$ there exists $\tau \in (0,1]$ such that for every $t \in [0,\tau],$ $X_t$ 
is transverse to the foliation of $\mathbb T^m.$
The remaining assertions follow from the theory 
developed in \cite{lawton}. In particular 
$\Lambda_t$ is irreducible because $X_t$ is connected. $\Lambda_t$ is
nontrivial because if it was an affine transformation of the product of
Bohr almost peridic sets on lower dimensional Euclidean spaces, then $X_t$ would either be disconnected or have dimension $> 1.$
\\ \\
We observe that $\varphi_t(\mathbb T) = \mathbb T^m \cap V$ where $V \subset \mathbb C^{*m}$ is the subvariety defined by the set of common zeros $z \in C^{*m}$ of Laurent polynomials:
\begin{equation}\label{variety} 
(z_j+t s_j)^{\gamma_k} (1+t s_k z_k)^{\gamma_j}
-(1+t s_j z_j)^{\gamma_k} (z_k+t s_k)^{\gamma_j}, \ \ j , k \in \{1,...,m\}, j \neq k.
\end{equation}
In general, if $n > 1,$ then $V$ is not defined by a subset consisting of $n$ of the $mn/2$ equations in (\ref{variety}), hence $V$ is not a complete intersection. Furthermore,
$\psi^{-1}(\varphi_t(\mathbb C^*)) \subset \mathbb R^n.$ This inclusion implies real-rootedness of related trigonometric maps.
\\ \\
A subset of $\mathbb R^n$ is called an infinite arithmetic progression if it equals $u + v\mathbb Z$ for $u \in \mathbb R^n$ and $v \in \mathbb R^n \backslash \{0\}.$
\begin{lem}\label{arith} 
If $t \in (0,\tau],$ then $\Lambda_t$ does not contain an infinite arithmetic progression. 
\end{lem}
Since $1,b_1,...,b_n$ are rationally independent $\psi := \rho_m \circ M$ is injective.
Assume to the contrary that $t \in (0,\tau]$ and 
$\Lambda_t$ contains an arithmetic progression 
$\Gamma_1.$ 
Then $\Gamma_1$ is contained in a finitely generated and thereore finite rank subgroup $\Gamma_2$ of $\mathbb R^n.$ Define 
$G_i := \psi(\Gamma_i), \, i = 1, 2.$
Then $G_2$ is a finite rank subgroup of $V.$ Therefore
the solution of Lang's conjecture (\cite{evertse}, Theorem 10.10.1) implies $G_2$ is 
contained in the finite union of translates of torus subgroups of $\varphi_t(\mathbb T)$ of dimension $< m.$
Since $s_1,...,s_m$ are distinct, $\varphi_t(\mathbb T)$ does not contain the translate of any torus subgroup of dimension $> 0,$ $G_2$ and hence $G_1$ is finite. Since 
$\psi$ is injective, $\Gamma_1$ is finite. This contradiction 
concludes the proof.
\begin{remark}\label{Lang}
Lang \cite{lang} conjectured that if $S$ is a finite rank subgroup of $C^{*m}$ 
and $V \subset C^{*m}$ is a proper variety, then $S \cap V$ is contained in a finite number of cosets
in $V$ of subgroups of $C^{*m}$ isomorphic to $C^{*k}$ where $k < m.$ 
Liardet \cite{liardet} proved this for $m = 2,$ and Laurent \cite{laurent} proved this for $m \geq 3.$
Ritt \cite{ritt} proved that if $p \in \mathcal T_1$ and 
equals $0$ on $u + v\mathbb Z,$ then there exists $p_1 \in \mathcal T_1$
such that $p(x) = 
p_1(x)\, \sin \pi (\frac{x}{v}-\frac{u}{v} ).$ This result implies that the set of real zeros of $p$ is a union $\Lambda_1 \cup \Lambda_2$ where $\Lambda_1$ does not contain an infinite arithmetic progression, and  $\Lambda_2 = \emptyset$ or $\Lambda_2$ is the union of a finite number of infinite arithmetic progressions.
\end{remark}
\begin{theo}\label{thm1}
$\Lambda_t$ is a Fourier quasicrystal for every $t \in [0,\tau].$
\end{theo}
Proof. The theory in \cite{lawton} implies that the Fourier-Bohr transform $\mathcal F_B(\Lambda_t)$ of $\Lambda_t$ equals zero outside the spectrum
$\Omega(\Lambda_t) \subset M^T\mathbb Z^m,$ and
\begin{equation}\label{FBmu1}
\mathcal F_B(\Lambda_t)(M^T\ell) = 
\int_{X_t)} z^{-\ell} d\theta, \ \ t \in [-,\tau], \, \ell \in \mathbb Z^m.
\end{equation}
where $X_t = \varphi_t(\mathbb T)$ is the $1$-cycle obtained by orienting $X_t$ by the counterclockwise orientation of $\mathbb T,$ and $d\theta$ is the differential $1$-form
\begin{equation}\label{dtheta}
d\theta = \frac{1}{2\pi i} \sum_{k = 1}^m
\alpha_k \frac{dz_k}{z_k}, \ \ \alpha := \frac{\Delta(\Lambda_0)}{[b^T \, 1]\ \gamma} \, \left[
\begin{array}{c}
b \\
1 
\end{array}
\right].
\end{equation}
Therefore $\alpha^T M = 0$ and $\alpha \cdot \gamma = \Delta(\Lambda_0).$ Since 
\begin{equation}\label{inequality}
|\mathcal F_B(\Lambda_t)(M^T\ell)| \leq 
\int_{\varphi_t(\mathbb T)} d\theta = \Delta(\Delta_0),
\ \ \ell \in \mathbb Z^m,
\end{equation}
it suffices to prove that there exists a polynomial $g$ such that
\begin{equation}\label{FBmu2}
|\, \Omega(\Lambda_t) \cap [-r,r]^n\, | \leq g(r), \  \ r > 0.
\end{equation}
Since $z_j(\zeta) = \frac{\zeta^{\gamma_j}-s_j}{1-s_j\zeta^{\gamma_j}}$ and
$\frac{dz_j}{d\zeta}(\zeta) = \gamma_j \zeta^{\gamma_j-1}
\frac{1-s_j^2}{1-s_j\zeta^{\gamma_j}}$ are
convergent power series in nonnegative powers of $\zeta,$
$\mathcal F_B(\Lambda_t)(M^T\ell) = 0$ whenever
all $\ell_j \leq -1.$ The fact that the Radon measure associated with $\Lambda_t$ is real-valued implies that
\begin{equation}\label{Hermitian}
\mathcal F_B(\Lambda_t)(-M^T\ell) = 
\overline {
\mathcal F_B(\Lambda_t)(M^T\ell)}, \ \ \ell \in \mathbb Z^m.
\end{equation}
Therefore  $\mathcal F_B(\Lambda_t)(M^T\ell) = 0$ whenever
all $\ell_j \geq 1.$
If $M^T\ell \in [-r,r]^n$ then
\begin{equation}\label{FBmu3}
\ell_j \in [\ell_mb_j - r,\ell_mb_j + r], \ \ j = 1,...,n.
\end{equation}
Let $\beta := \max \{1/b_j : j = 1,...,n\}.$ 
If $|\ell_m| \geq \max \{1, \beta(r+1)\}$
then either all $\ell_j \geq -1$ or all $\ell_j \geq 1$
and $\mathcal F_B(\mu)(M^T\ell) = 0.$
Then 
\begin{equation}\label{g}
g(r)) := 2(\beta r + 1 + \beta)(2r+1)^n
\end{equation}
satisfies inequality (\ref{FBmu2}) and the proof is concluded.
\begin{remark}
For $n = 1, s_1 = -\frac{1}{3}, s_2 = 0, \gamma_1 = 2, \gamma_2 = 1$
\begin{equation}\label{KS}
q_1(z_1,z_2) = 
1 - \frac{1}{3}z_1 + \frac{1}{3}z_2^2- z_1z_2^2
\end{equation}
is equivalent (under the automorphism 
$[z_1, z_2]^T \mapsto [z_1, z_2^{-1}]^T$) 
to the Laurent polynomial that Kurasov and Sarnak used (\cite{kurasovsarnak}, Equation 41) 
to construct the first nontrivial FQ on $\mathbb R.$
\end{remark}
%
%
\section{Amoebas and Real-Rootedness}\label{sec4}
\subsection{Background}\label{subsec4.1}
In 1952 in their study of phase transitions in statistical physics Lee and Yang introduced a family of $q \in \mathcal L_m$ and proved (\cite{leeyang}, Appendix II) that
\begin{equation}\label{leeyangcond}
	q([z_1,...,z_n]^T) \neq 0 \  \text{whenever all} \ |z_j| < 1 \ \text{or all}  \ |z_j| > 1,
\end{equation}
so $c(z) := q([z,...,z]^T)$ has all its zeros in the unit circle. This family is 
\begin{equation}\label{leeyangpoly}
	q([z_1,...,z_n]^T) := \sum_J \prod_{j \in J} \left( z_j \prod_{k \notin J} a_{j,k} \right)
\end{equation}
where $A := (a_{j,k}) \in \mathbb C^{n \times n}$ is symmetric with $0 < |a_{j,k}| < 1$ and $J$ is summed over all subsets of $\{1,...,n\}.$ 
Subsequently more general polynomials satisfying (\ref{leeyangcond}) were studied and named Lee-Yang polynomials. 
In 2010 Ruelle completely characterized multiaffine Lee-Yang polynomials (degree $1$ in each variable) (\cite{ruelle}, Theorem 3), higher degree in each variable Lee-Yang polynomials (\cite{ruelle}, Section 10), and interpreted the latter as temperature dependent partition functions involving highr spins
(\cite{ruelle}, Section 11).
\\ \\
In 2005 Passare and Tsikh (\cite{passaretsikh}, Theorem 2) interpreted (\ref{leeyangcond}) in terms of the amoeba concept to prove that polynomials defined by (\ref{leeyangpoly}) satisfy (\ref{leeyangcond}). To explain this concept, define the logarithmic moment map 
$Log : \mathbb C^{*m} \mapsto \mathbb R^m$
\begin{equation}\label{Log}
	Log \left(\, [z_1,...,z_m]^T\, \right) := [\, \ln |z_1|,...,\ln |z_m|\, ]^T.
\end{equation}
$Log$ is a homomorphism from $C^{*m}$ onto 
$\mathbb R^m$ with kernel $\mathbb T^m.$ 
A subset $V$ of $\mathbb C^{*m}$ is called an algebraic variety if it is the set of common zeros of a set of Laurent polynomials. By the Noetherian property of the ring of Laurent polynomials, this set of polynomials can be chosen to be finite. $V$ is called a hypersurface if the set has one element and then $V$ has codimension $1.$ $V$ is called a complete intersection if it equals the intersection of $k$ hypersurfaces and has codimension $k.$
The amoeba of an algebraic variety $V \subset \mathbb C^{*m}$ is
\begin{equation}\label{amoeba1}
 \mathcal A(V) := Log (V).
\end{equation}
Gelfand, Kapranov and Zelevinsky (\cite{gelfand}, p. 194, Definition 1.3) defined amoebas for hypersurfaces and proved (\cite{gelfand}, p. 195, Corollary 1.6) that the complement of the amoeba of a hypersurface is a finite union of open convex subsets. Several researchers \cite{forsberg,forsbergpassaretsikh,passaretrullgard} extended their work.. 
Clearly $q \in \mathcal L_m$ is  a Lee-Yang polynomial iff
$\mathcal A(\Lambda(q)) \cap \pm \mathbb R_+^m = \emptyset,$ or equivalently
\begin{equation}\label{amoebaeqn2}
\mathcal A(\Lambda(q)) \cap M \mathbb R = \{0\}, \ \ M \in \mathbb R_+^m.
\end{equation}
In \cite{alon1} Alon, Cohen and Vinzant used amoebas to prove that every real-rooted  $p \in \mathcal T_1$ admits the representation
\begin{equation}\label{representation1}
p  = q \circ \rho_m \circ M 
\end{equation}
where $m = \text{dim}_{\mathbb Z}(\Omega(p)),$ 
$q \in \mathcal L_m$ is a Lee-Yang polynomial, and
$M \in R_+^m.$ We observe that (\ref{amoebaeqn2}) holds if $M$ is replaced by a matric $M_1$ sufficiently close to $M,$ a condition we call $M$-stable.
\\ \\
This section proves Theorem \ref{thm2}, which represents uniformly generic real-rooted trigonometric maps $P \in \mathcal T_n^n$ by $M$-stable Laurent maps. Here $M \in \mathbb R^{m \times n}.$ 
This concept generalizes the concept of Lee-Yang polynomials.
\subsection{$M$-stable Laurent Maps}\label{subsec4.2}
\begin{defi}\label{stable1}
For a rank $n$ matrix $M \in \mathbb R^{m \times n},$ $Q \in \mathcal L_m^n$ is $M$-stable if
\begin{equation}\label{stable2}
	\mathcal A(\Lambda(Q)) \cap M^{\prime} \mathbb R^n = \{0\}
\end{equation}
for all $M^{\prime}$ in some open neighborhood of $M.$ 
\end{defi}
\begin{lem}\label{lemma13} 
If $M_1, M_2 \in \mathbb R^{m \times n}$ have rank $n$ and $Q_1$ is $M_1$-stable, then $Q_1 \sim Q_2$ for some $M_2$ stable $Q_2.$ Here the equivalence relation $\sim$ was defined in the first paragraph in Subsection \ref{subsec2.6}.
For $n = 1$ we can choose $M_2 \in \mathbb R_+^m$ and $q_2$ to be a Lee-Yang Laurent polynomial.
\end{lem}
Proof. $GL(m,\mathbb Z)$ acts on the Grassmann manifold $Gr(n,m)$ consisting of $n$-dimensional subspaces of $\mathbb R^m.$ An element
in $Gr(n,m)$ is rational if it has the form $N\mathbb R^n$ where $N \in \mathbb Z^{m \times n}$ has rank $n.$
The existence of the Smith normal form \cite{newman, smith} for integer matrices implies that $GL(m,\mathbb Z)$ acts transitively on the set of rational elements in $Gr(n,m).$ Since the set of rational elements is dense, and the compact group $O(m,\mathbb R)$ acts transitively on $Gr(m,mn),$ every orbit in $Gr(n,m)$ under the action of $GL(m,\mathbb Z)$ is dense. The hypothesis implies that there exists an open $O \subset Gr(n,m)$ such that 
$M_1\mathbb R^n \in O$ and
\begin{equation}\label{Q1}
	\mathcal A(\Lambda(Q_1)) \cap S = \{0\}, \ \ S \in O.
\end{equation}
Since every orbit is dense, there exists $U \in GL(m,\mathbb Z)$ such that $UM_2\mathbb R^n \in O$
so $U^{-1}O \subset Gr(n,m)$ is open and $M_2\mathbb R^n \in U^{-1}O.$ Define $Q_2 := Q_1 \circ \widetilde U.$ Then
\begin{equation}\label{Q2}
	\mathcal A(\Lambda(Q_2)) \cap S = \{0\}, \ \ S \in U^{-1}O.
\end{equation}
Since $M_2^{\prime} \mathbb R^n \in U^{-1}O$ for all $M_2^{\prime}$ in some open neighborhood of $M_1,$ $Q_2$ is $M_2$-stable. This proves the first assertion. 
To prove the second assertion assume that $M_1 \in \mathbb R^m$ is nonzero and that a Laurent polynomial $q_1$ is $M_1$-stable. Then there exist
an open $O \subset Gr(1,m)$ and $u \in \mathbb Z^m$ such that $gcd(u_1,...,u_m) = 1,$
$u \mathbb R \in O,$
$M_1\mathbb R \in O,$ and
\begin{equation}\label{q1}
	\mathcal A(\Lambda(q_1)) \cap S = \{0\}, \ \ S \in O.
\end{equation}
The Smith normal form implies that
there exists $U_1 \in GL(m,\mathbb Z)$ whose first column equals $u.$ By adding sufficiently large positive integer multiples of 
$u$ to the last $m-1$ columns of $U_1$ we construct $U \in GL(m,\mathbb Z)$ such that 
\begin{equation}\label{v}
	Uv\mathbb R \in O, \ \ v \in (\{0\} \cup \mathbb R_+)^m.
\end{equation}
Define $q_2 := q_1 \circ \widetilde U.$ Then for every $v \in (\{0\} \cup \mathbb R_+)^m,$
\begin{equation}\label{q2}
	\mathcal A(\Lambda(q_2)) \cap v \mathbb R = U^{-1}[\mathcal A(\Lambda(q_1)) \cap Uv \mathbb R]
= U^{-1}\{0\} = \{0\}
\end{equation}
hence $q_2$ is a Lee-Yang Laurent polynomial.
\subsection{Proof of Proposition 7}\label{subsec4.3}
The following result is an extension to Laurent maps of the result for Laurent polynomials by Alon, Cohen and Vinzant (\cite{alon1}, Corollary 3.3).
%
%
\begin{prop}\label{prop7}
If $Q \in \mathcal L_m^n$ and $M \in \mathbb R^{m \times n}$ has rank $n$ and $P := Q \circ \rho_n \circ M$ is uniformly generic, then
\begin{equation}\label{ImP}
 \mathcal A(\Lambda(Q)) \cap M\mathbb R^n = -2\pi M \, \overline {\Im \Lambda(P)}.
\end{equation}
\end{prop}
Proof. Define
\begin{equation}\label{GM}
G(M) := \overline {\rho_m(M\mathbb C^n)}.
\end{equation}
Since $G(M) = \{z \in \mathbb C^{*m} : \ell^T M = 0 \implies z^\ell - 1 = 0 \}$
is a variety and 
$\mathcal A(\Lambda(Q)) \cap M\mathbb R^n 
= \mathcal A(\Lambda(Q) \cap G(M)) \cap M\mathbb R^n,$ it suffices to prove that
\begin{equation}\label{ImPext}
 \mathcal A(\Lambda(Q) \cap G(M)) \cap M\mathbb R^n = -2\pi M \, \overline {\Im \Lambda(P)}.
\end{equation}
If $\lambda = [\lambda_1,...,\lambda_n]^T \in \Lambda(P),$
then $\rho_m(M\lambda) \in \Lambda(Q) \cap G(M)$ hence 
$$
Log (\rho_m(M \lambda)) = -2\pi M \, \Im \lambda \in \mathcal A(\Lambda(Q) \cap G(M)) \cap M\mathbb R^n,
$$
so 
$$-2\pi M \Im \Lambda(P) \subset \mathcal A(\Lambda(Q))  \cap M\mathbb R^n.$$ 
Since $\mathcal A(\Lambda(Q)) \cap M\mathbb R^n$ is closed, 
$$-2\pi M \overline {\Im \Lambda(P)} \subset \mathcal A(\Lambda(Q)) \cap M\mathbb R^n.$$
Proving the reverse inclusion is more challenging. 
We use proof by contradiction. Assume that there exists
\begin{equation}\label{u}
 u \in \mathcal A(\Lambda(Q) \cap G(M)) \cap M\mathbb R^n \, \backslash
-2\pi M \overline {\Im \Lambda(P)}.
\end{equation}
Since $-2\pi M \overline {\Im \Lambda(P)}$ is closed 
there exists $\delta > 0$ such that 
\begin{equation}\label{empty}
B_m(u,2\delta) \cap -2\pi M \, \Im \Lambda(P) = \emptyset
\end{equation}
therefore $P$ does not vanish on the set 
\begin{equation}\label{H}
H_{2\delta} := \{ z \in \mathbb C^n : -2\pi M \Im z \in B_m(u,2\delta)\, \}.
\end{equation}
Proposition \ref{prop4} implies rank $M = n$ hence $\Im H_{2\delta}$ is bounded. Since $H_{2\delta}$ is
invariant under translation by elements in 
$\mathbb R^n,$ $H_{2\delta}$ is a tube domain and Corollary \ref{corProp3} implies that $P \in B(H_{2\delta}).$ Moreover the Hausdorff distance
\begin{equation}\label{distance}
 d(\overline {H_\delta}, \mathbb C^n \backslash H_{2\delta}) > 0.
\end{equation}
Proposition \ref{prop5} implies there exists $\epsilon > 0$ such that
\begin{equation}\label{Pbound}
 |P(z)| \geq \epsilon, \ \ z \in H_\delta.
\end{equation}
(\ref{u}) and the fact that $\rho_m : \mathbb C^m  \mapsto \mathbb C^{*m}$ is surjective implies there exists
$w \in \mathbb C^m$ such that 
\begin{equation}\label{w}
 \rho_m(w) \in \Lambda(Q) \cap G(M)
\end{equation}
\begin{equation}\label{uw}
 u = Log \, \rho_m(w) = -2\pi \Im w.
\end{equation}
(\ref{u}) implies there exists  $v \in \mathbb R^n$ such that
\begin{equation}\label{v}
 u = Mv,
\end{equation}
and (\ref{uw}) and (\ref{v}) implies $\rho_m(i \Im w) \in G(M)$ therefore
(\ref{w}) implies
\begin{equation}\label{Rew}
 \rho_m(\Re w) = \rho_m(w - i \Im w) = \rho_m(w) \rho_m(- i \Im w) \in G(M).
\end{equation}
For $\tau \in \mathbb R^n$ define $P_\tau \in \mathcal T_n^n$ by
\begin{equation}\label{Ptau}
 P_\tau(z) := Q \circ \rho_m((\Re w - M\tau) + Mz), \ \ z \in \mathbb C^n.
\end{equation}
(\ref{w}), (\ref{uw}), (\ref{v}) and (\ref{Ptau}) imply
\begin{equation}\label{root}
 P_\tau \left(\tau + \frac{1}{2\pi i}v \right) = 
Q \circ \rho_m \left(\Re w +\frac{1}{2\pi i}u) \right) = 
Q \circ \rho_m(\Re w +i \Im w) = 0.
\end{equation}
(\ref{GM}) and (\ref{Rew}) imply
there exists a sequence $\tau_k \in \mathbb R^n$ such that
\begin{equation}\label{Ptaumult}
 \lim_{k \rightarrow \infty} \rho_m(\Re w - M\tau_k) = 1 \in \mathbb T^m.
\end{equation}
Therefore
\begin{equation}\label{limP}
 \lim_{k \rightarrow \infty} P_{\tau_k}(z) = Q(\rho_m(Mz)) = P(z)
\end{equation}
and the convergence is uniform over $H_\delta.$ Choose
$k$ sufficiently large that 
\begin{equation}\label{approx}
 |P_{\tau_k}(z) -  P(z)| \leq \frac{\epsilon}{2}, \ \ z \in H_\delta.
\end{equation}
(\ref{Pbound}) and (\ref{approx}) imply
\begin{equation}\label{Ptaukbound}
 |P_{\tau_k}(z)| \geq \frac{\epsilon}{2}, \ \ z \in H_\delta
\end{equation}
so $P_{\tau_k}$ has no roots in $H_\delta.$ This contradicts (\ref{root})
and concludes the proof.
\subsection{Proof of  Theorem 2}\label{subsec4.4}
\begin{theo}\label{thm2}
If $P \in \mathcal T_n^n$ is uniformly generic and real-rooted and $m = \text{dim}_{\ \mathbb Z} \Omega(P),$ then there exists a rank $n$ matrix $M \in \mathbb R^{m \times n}$ with $\text{dim}_{\mathbb Z} M = m$ and an $M$-stable
$Q \in \mathcal L_m^n$ with
\begin{equation}\label{representation2}
P = Q \circ \rho_m \circ M.
\end{equation}
\end{theo}
\begin{remark}\label{unidimension}
In \cite{kurasovsarnak} Kurasov and Sarnak observed that if $Q \in \mathcal L_m$
satisfies $\pm \mathbb R_+^m \subset \mathbb R^m \backslash \mathcal A(\Lambda(Q))$ and
$M \in \mathbb R_+^m$ then $P := Q \circ \rho_m \circ M \in \mathcal T_m$ is real-rooted. 
In (\cite{alon1},Theorem 1.1) the authors proved that every real-rooted trigonometric polynomial admits such a representation. If $M_1$ is sufficiently close to $M,$ then
$M_1 \in \mathbb R_+^m$ hence $P_1 := Q \circ \rho_m \circ M_1$ is real-rooted and $P$ is $M$-stable,
\end{remark}
The proof of Theorem \ref{thm2} requires the following result.
\begin{lem}\label{localstab} 
If $\mathcal A(Q) \cap M\mathbb R^n  = \{0\}$ and
rank $M = n,$ then
there exists $\epsilon > 0$ such that
\begin{equation}\label{zero2}
\mathcal A(\Lambda(Q)) \cap M_1\mathbb R^n \cap B_n(0,\epsilon) = \{0\}
\end{equation}
for $M_1$ sufficiently close to $M.$
\end{lem}
Proof. Let $Q = [q_1,...,q_n]^T$ and let $S_m$ be the boundary of $B_n(0,1) \subset \mathbb R^m$
and $S := S_m \cap M\mathbb R^n.$ The hypothesis implies that for every 
$s \in S$ the open ray 
\begin{equation}\label{ray}
s\mathbb R_+ \subset \mathbb R^m \, \backslash  \mathcal A(\Lambda(Q)) = \bigcup_{j=1}^n \left[ \, \mathbb R^m \, \backslash  \mathcal A(\Lambda(q_j)) \, \right].
\end{equation}
Since each $\mathbb R^m \, \backslash  \mathcal A(\Lambda(q_j))$
is a finite union of open convex subsets there must exist
$j \in \{1,...,n\}$ and $\epsilon(s) > 0$ such that the line segment
$s (0,\epsilon(s)]$ is contained in one of these open convex subsets. 
Therefore there exists an open cone
$C_s \subset \mathbb R^m$ such that 
$$
 s\mathbb R_+ \subset C_s
$$
and
$$C_s \cap B_n(0,\epsilon(s)) \subset \mathbb R^m \, \backslash  \mathcal A(\Lambda(Q)).$$
Since $\{C_s, s \in S\}$ is an open cover of $S,$ there exists
a finite collection $C_{s_1},...,C_{s_L}$ whose union contains $S.$
Define $\epsilon := \min \{\epsilon_{s_j}, j = 1,...,L\}$ and
$O := \bigcup_{j=1}^L C_s.$
Then
$$
O \cap B_m(0,\epsilon) \subset \mathbb R^m \, \backslash  \mathcal A(\Lambda(Q)),
$$
hence if $M_1$ is sufficiently close to $M$ then
$$
		M_1 \mathbb R^n \subset O \cup \{0\}.
$$
The proof of Lemma \ref{localstab}  is concluded since
$$
\mathcal A(\Lambda(Q)) \cap M_1\mathbb R^n \cap B_m(0,\epsilon) 
\subset \{0\} \cup \left[ \mathcal A(\Lambda(Q)) \cap O \cap B_m(0,\epsilon) \right] = \{0\}.
$$
We complete the proof of Theorem \ref{thm2} by contradiction. 
Assume that $P$ is uniformly generic and real-rooted and let $P = Q \circ \rho_n \circ M$ be the representation in Proposition \ref{prop3}. Since $P$ is generic, rank $M = n.$
Assume to the contrary that $P$ is not $M$-stable. Then Lemma \ref{localstab}
implies that for every neigborhood $U_M$ of $M$ there exists $M_1 \in U_M$
and $u \in \mathcal A(\Lambda(Q)) \cap M_1\mathbb R^n \backslash B_m(0,\epsilon).$ 
We may assume that is $U_M$ sufficiently small $P_1 := Q \circ \rho_m \circ M_1$ is uniformly generic 
so there exists $R > 0$ with
\begin{equation}\label{imrootbound}
	-2 \pi M_1 \overline {\Im \Lambda(P_1)} \subset B_m(0,R). 
\end{equation}
Proposition \ref{prop7} then implies
\begin{equation}\label{ubound}
	\epsilon < |u| < R.
\end{equation}
Let $S_{\epsilon}$ and $S_{R}$ be spheres 
radii $\epsilon$ and $R$ in $M\mathbb R^n$ 
oriented in a positive direction to define cycles in
$M_1\mathbb R^n \backslash \{0\}$ and in $\mathbb R^m \backslash \mathcal A(\Lambda(Q)).$
Then
the homology classes $H_{n-1}(S_{\epsilon}) \neq H_{n-1}(S_R)$ in 
$H_1(M_1\mathbb R^n \backslash \mathcal A(\Lambda(Q))$ because the point $u$ lies between
the spheres.
However, these cycles are homotopic in $\mathbb R^m \backslash \mathcal A(\Lambda(Q))$ as the following argument shows. Define 
\begin{equation}\label{Mt}
	((1-t)M_1+tM) \in \mathbb R^{m \times n}, \ \ t \in [0,1].
\end{equation}
Then $M_tS_{\epsilon}$ and $M_tS_{R}$
define homotopies that deforms the $(n-1)$-cycles
$S_{\epsilon}$ and $S_R$ to homologous cycles in 
$M\mathbb R^n \backslash \{0\}.$ 
Therefore the cycles $S_{\epsilon/2}$ and $S_{R}$ are homologous in $\mathbb R^m \backslash \mathcal A(\Lambda(Q))$
so $H_{n-1}(S_{\epsilon}) = H_{n-1}(S_R)$ in $\mathbb R^m \backslash \mathcal A(\Lambda(Q)).$ 
This means that the homology inclusion 
\begin{equation}\label{hominc}
H_{n-1}(M_1\mathbb R^m \backslash \mathcal A(\Lambda(Q))) \mapsto 
H_{n-1}(M\mathbb R^m \backslash \mathcal A(\Lambda(Q)))
\end{equation}
is not injective.
This violates (\cite{bushuevatsikh}, Theorem 1) and concludes the proof. 
\begin{remark}\label{conj}
The result of Bushueva and Tsikh (\cite{bushuevatsikh}, Theorem 1) was conjectured and proved under additional hypotheses by Henriques \cite {henriques}.
\end{remark}
%
%
\section{Sufficient Conditions for Fourier Quasicrystals on $\mathbb R^n$}\label{sec5}
This section proves the following result:
\begin{theo}\label{thm3}
If $P \in \mathcal T_n^n$ is real-rooted and uniformly generic then $div\, P$ is a Fourier Quasicrystal on $\mathbb R^n.$
\end{theo}
\subsection{Approximation Strategy}\label{subsec5.1}
Since $\mathcal S_c(\mathbb R^n)$ is dense in 
$\mathcal S(\mathbb R^n)$ it suffices to prove the existence of
a discrete $S \subset \mathbb R^n$ and 
$a_s \in \mathbb C^*, s \in S,$ such that both
$\zeta := \sum_{s \in S} a_s \delta_s$ and
$|\zeta|$ are tempered Radon measures and
\begin{equation}\label{step1}
<\widehat {div P},h> \, = \, <\zeta, h>, \ \ h \in \mathcal S_c(\mathbb R^n).
\end{equation}
Assume that $P \in \mathcal T_n^n$ is real-rooted and generic and $m := \text{dim}_{\mathbb Z} \Omega(P).$ Theorem \ref{thm2} implies there exists a rank $n$ matrix
$M \in \mathbb R^{m \times n}$ and an $M$-stable Laurent map $Q \in \mathcal L_m^n$ such that
$P = Q \circ \rho_m \circ M.$ Define a sequence 
$P_k \in \mathcal T_n^n$ 
\begin{equation}\label{Pk1}
	P_k := Q \circ \rho_m \circ M_k.
\end{equation}
where $M_k \in \mathbb Q^{m \times n}$ 
converges to $M.$  
Since $Q$ is $M$-stable we may assume that each $P_k$ is real-rooted and generic.
\begin{lem}\label{densconv} 
	$\Delta(div\, P) = \lim_{k \rightarrow \infty} \Delta(div\, P_k)$
\end{lem}
Proof. Since each component of $\mathcal N(P)$ (resp. $\mathcal N(P_k))$ equals $M^T$ (resp. $M_k^T$) 
times the corresponding component of $\mathcal N(Q),$ 
it follows that 
$$V(\mathcal N(P)) = \lim_{k \rightarrow \infty}
V(\mathcal N(P_k))$$ 
where convergence is componentwise with respect to the Hausdorff metric. The conclusion follows (\ref{dendivP1}).
%
%
\begin{prop}\label{prop8} 
$<\widehat {div\, P}, h> \, = \, \lim_{k \rightarrow \infty} \, <\widehat  {div\, P_k}, h>, \ \ h \in \mathcal S_c(\mathbb R^n).$
\end{prop}
Proof. For $r > 0,$ let $f_r$ be the characteristic function of $B_n(0,r).$ (\ref{duality2}) implies
$$| <\widehat {div\, P} - \widehat {div\, P_k}, h> | \leq | <f_r \, div\, P - f_r \, div\, P_k, \widehat h> |$$
$$ + \ | <(1-f_r) div\, P, \widehat h> |+ | <(1-f_r) div\, P_k, \widehat h> |.$$
Let $T \subset \mathbb C^n$ be a tube domain containing $\mathbb R^n.$ 
$P$ and $P_k$ are Bohr almost periodic in $T$ and $P_k$ converges to
$P$ uniformly on compact subsets of $T.$  Lemma \ref{rouche} implies that 
$f_r \, div \, P_k$ converges $f_r \, div \, P.$ Therefore
the first term on the right above converges to $0.$  $\Delta(div\, P) < \infty$
implies $div\, P(B_n(0,s)) = O(s^n).$ Therefore
\begin{equation}\label{split2}
|\, < (1-f_r) div\, P, \widehat h> \, | \leq 
\int_{s=r}^\infty \sup_{|x| = s} |\widehat h(x)| \, d (div\, P(B_n(0,s)))
\end{equation}
converges to $0$ as $r \rightarrow \infty$ since $\sup_{|x| = s} |\widehat h(x)|$ has fast decay. Similarly the third term on the right above converges to $0$ 
as $r \rightarrow \infty.$
\\ \\
Express 
$M_k = \beta_k N_k,$ 
where 
$\beta_k^{-1} \in \mathbb Z_+$ is the least common multiple of the denominators of the entries in $M_k$ 
and 
$N_k \in \mathbb Z^{m \times n}.$ 
(\ref{rhomtildeN}) gives
\begin{equation}\label{Pk2}
	P_k = Q_k \circ \rho_n \circ \beta_k \, I_n
\end{equation}
where
\begin{equation}\label{Qk}
	Q_k = Q \circ \widetilde {N_k^T} \in \mathcal L_n^n. 
\end{equation}
Since $P_k$ and hence $Q_k$ is generic Lemma 
\ref{bernshtein}
implies $\Lambda(Q_k)$ is finite. Then since $P_k$ is real-rooted
$\Lambda(Q_k) \subset \mathbb T^n$
so there exists a finite $X_k \subset \mathbb R^n$ with
\begin{equation}\label{rootsQk4}
	\Lambda(Q_k) = \{\rho_n(x) : x \in X_k\}
\end{equation}
For $x \in X_k$ define
\begin{equation}\label{Lsigma}
	L_x := \beta_k^{-1}\mathbb Z^n + \beta_k^{-1}x
\end{equation}
If $y \in L_x$ then $P_k(y) = 0$ and the multiplicities satisfy
\begin{equation}\label{mult}
	m_y(P_k) = m_{\rho_n(x)}(Q_k)
\end{equation}
since $\rho_n \circ \beta_kI_n(y) = \rho_n(x)$ and multiplicity is preserved under the holomorphic local homeomorphisms.
Therefore
\begin{equation}\label{divPk}
	div P_k = \sum_{x \in X_k} m_{\rho_n(x)}(Q_k)
\sum_{y \in L_x} \delta_y.
\end{equation}
The Poisson summation formula Lemma \ref{PSF1}
and its consequent (\ref{PSF3}) imply
\begin{equation}\label{FTdivPk}
	\widehat {div P_k} = 
	\beta_k^n \sum_{\lambda \in \lambda(Q_k)} m_{\Lambda}(Q_k)
\sum_{\ell \in \mathbb Z^n} \lambda^{-\ell} \delta_{\beta_k \ell}.
\end{equation}
\begin{defi}\label{Q1}
	For $h \in \mathcal S_c(\mathbb R^n)$ and $k \in \mathbb Z_+$ define $R_{h,k} \in \mathcal L_n$ by
\begin{equation}\label{Rhk}
R_{h,k}(z) :=  \sum_{\ell \in \mathbb Z^n} h(\beta_k\ell)\, z^{-\ell}, \ \ z \in \mathbb C^{*n}.
\end{equation}
\end{defi}
(\ref{FTdivPk}) and (\ref{Rhk}) imply
\begin{equation}\label{FTdivPkh}
<\widehat {div P_k}, h> \, = \, 
\beta_k^n \sum_{\lambda \in \Lambda(Q_k)} m_\lambda(Q_k)\, R_{h,k}(\lambda), \ \ h \in  \mathcal S_c(\mathbb R^n). 
\end{equation}
Define differential $n$-forms on $\mathbb C^{*n} \backslash \Lambda(q_k)$
\begin{equation}\label{omega}
	\omega_{h,k} := \left(\frac{\beta_k}{2\pi i}\right)^n R_{h,k} 
	\frac{dq_{k,1}}{q_{k,1}} \wedge \cdots \wedge \frac{dq_{k,n}}{q_{k,n}}, \ \ h \in  \mathcal S_c(\mathbb R^n), k \in \mathbb Z_+
\end{equation}
where $Q_k = [q_{k,1},...,q_{k,n}]^T$ and $q_ k := q_{k,1} \cdots q_{k,n}.$
\\ \\
Define $n$-chains in $\mathbb C^{*n} \backslash \Lambda(q_k)$
\begin{equation}\label{omega}
	\Gamma_k = \sum_{\lambda \in \Lambda(Q_k)} \Gamma_\lambda, \ \ k \in \mathbb Z_+
\end{equation}
where $\Gamma_\lambda$ is the Grothendieck cycle defined in (\ref{Gcycle})
with $f_j$ replaced by $q_{k,j}.$
\begin{lem}\label{intG} 
\begin{equation}\label{GrowInt}
<\widehat {div P_k}, h>\, = \int_{\Gamma_k} \omega_{h,k} 
\end{equation}
\end{lem}
Proof. Follows from (\ref{FTdivPkh}) and (\ref{growint}) by letting $\epsilon_j \rightarrow 0$ in (\ref{Gcycle}).
\subsection{Grothendieck Residues and Toric Geometry}\label{subsec5.2}
In (\cite{gelfondkhovanskii2}, 1.2) Gelfond and Khovanskii introduced the following concept:
\begin{defi}
Let $[K_1,...,K_n]^T$ be an $n$-tuple of  polytopes in $R^n.$ Then each face $F$ of their Minkowski sum $K := K_1 + \cdots + K_n$ is a Minkowski sum $F = F_1 + \cdots  + F_n$ where each $F_j$ is a unique face of $K_j$ called a summand of $F.$ The tuple is called unfolded if every face of $K$ has at least one vertex among it summands. Otherwise it is called folded.
\end{defi}
The Laurent maps $Q \in \mathcal L_2^2$ in
in \ref{subsec2.6}
have Newton polytope tuples
$\mathcal N(Q) = [K_1, K_2]^T$ where $K_1$ and $K_2$ are identical squares. Therefore $K_1+K_2$ is a square each of whose sides has summands consisting of sides not vertices. Therefore $[K_1, K_2]^T$ is folded.
We record the following observation without proof. 
$[K_1,...,K_n]^T$ is unfolded iff for every $u \in \mathbb R^n$ there exists $j \in \{1,...,n\}$ and a vertex $v \in K_j$
such that
\begin{equation}\label{unfolded}
 u \cdot v < u \cdot w, \ \ \text{for every vertex } w \in K_j \backslash \{v\}.
\end{equation}
It is easy to show that if $Q \in \mathcal L_n^n$ and the components of
$\mathcal N(Q)$ satisfy (\ref{unfolded}) then $Q$ is generic.
\\ \\
If a tuple of polytopes is unfolded then to each vertex $v$ of their sum $K$ is associated an integer $C(v)$ called its combinatorial coefficient, see \cite{gelfond2,gelfondkhovanskii1,gelfondkhovanskii2}.
If $P : \mathbb C^n \mapsto \mathbb C^n$ is a trigonometric map and $\mathcal N(P)$ is unfolded, then for $P_1$ sufficiently close to $P$ (with respect to its frequencies and corresponding amplitudes) $\mathcal N(P_1)$ is unfolded 
and $C(v_1) = C(v)$ whenever $v_1$ is a vertex of the sum of the polytopes in
$\mathcal N(P_1),$ and $v$ is the vertex of the sum of the polytopes in
$\mathcal N(P)$ that is closest to $v_1.$
We may assume that each $\mathcal N(P_k)$ is unfolded so $\mathcal N(Q_k)$ is unfolded since its entries differ from those in $\mathcal N(P_k)$ by a constant multiple.
\\ \\
For every $r \in \mathbb R_+^n$ define the Newton cycle $N_r := r\mathbb T^n$ oriented so 
$d(arg\, z_1) \wedge \cdots \wedge d(arg \, z_n)$ is positive.
%
%
\begin{prop}\label{prop9}
If $\mathcal N(Q_k)$ is unfolded and $q_k = q_{k,1} \cdots q_{k,n},$ 
then for every $v \in \mathcal V(q_k)$ there exists $r(k,v) \in \mathcal R_+^n$ such that 
\begin{enumerate}
\item $N_{r(k,v)} \subset \mathbb C^{*n} \backslash \Lambda(q_k),$
\item In the restriction of $q_k$ to $N_{r(k,v)}$ the monomial corresponding to $v$ 
dominates the other monomials in $q_k$ so that for some $q_{k,v} \in \mathbb C^*,$
\begin{equation}
q_k(z) = q_{k,v}z^v \, (1 - g_{k,v}(z))
\end{equation} 
such that $|g_{k,v}(z)| \leq \frac{1}{2}$ for all $z \in N_{z(v,k)}.$
\item $\Gamma_k$ is homologous to the Newton $n$-chain defined by
\begin{equation}
N_k := (-1)^n \sum_{v \in \mathcal V(q_k)} C(k,v) \, N_{r(k,v)},
\end{equation}
where $C(k,v)$ is the combinatorial coefficient for the vertex $v.$ 
\item Then Lemma (\ref{intG}) and Stokes theorem imply
\begin{equation}\label{stokes}
<\widehat {div\, P_k}, h> \, = \int_{N_k} \omega_{h,k}.
\end{equation}
\end{enumerate}
\end{prop}
Proof. These results are derived in (\cite{gelfondkhovanskii2}. 
\subsection{Proof of Proposition 10}\label{subsec5.3}
%
%
\begin{prop}\label{prop10}
If $P \in \mathcal T_n^n$ is real-rooted and $\mathcal N(P)$ is unfolded, then $div\, P$ is a Fourier Quasicrystal on $\mathbb R^n.$
\end{prop}
The proof of this will occupy the remainder of this subsection.
Define $H_{h,k}(z) := R_{h,k}(z)\, \frac{\beta_k^n z_1 \cdots z_n \det \left(\partial_{z_i} q_{k,j}\right)}{q_k}.$
Then $H_{h,k}$ is a rational function holomorphic in $\mathbb C^{*n} \backslash \Lambda(q_k).$ A computation gives
\begin{equation}\label{omegaH}
	 \omega_{h,k} = H_{h,k} \left(\frac{1}{2\pi i}\right)^n \frac{dz_1}{z_1} \wedge \cdots \wedge \frac{dz_n}{z_n}
\end{equation}
therefore
\begin{equation}\label{Hhkv}
	 \int_{N_{r(k,v)}} \omega_{h,k} = H_{h,k,v}
\end{equation}
where $H_{h,k,v}$ is the constant term in the Laurent expansion of the restriction
\begin{equation}\label{LaurentH}
	 H_{h,k}(z)|_{N_{r(k,v)}} = 
\frac{R_{h,k}(z) \, \beta_k^n \, z_1 \cdots z_n \det\left(\partial_{z_i} q_{k,j}(z)\right)}
{q_{k,v}z^v}\, 
\sum_{m = 0}^\infty g_{k,v}(z)^m.
\end{equation}
Express (\ref{LaurentH}) in terms of the components of $P_k,$ $p_{k,j}(w) = q_{k,j}(z)$ where $z = \rho_n(\beta_k w)$ for $w \in \mathbb C^n.$ The chain rule gives the trigonometric series
\begin{equation}\label{wedgep}
	H_{h,k}(\rho_n(\beta_kw)) = \frac{R_{h,k}(\rho_n(\beta_kw)) \, \det \left(\partial_{w_i} p_{k,j}(w)\right)}{q_{k,v}e^{2\pi i \beta_k v \cdot w}}\sum_{m = 0}^\infty g_{k,v}(\rho_n(\beta_kw))^m.
\end{equation}
when $\rho_n(\beta_k w) \in N_{r(k,v)},$ or equivalently $\ln |r(k,v)| = -2\pi\Im w.$ $H_{h,k,v}$ is the constant term in this series.
\begin{lem}\label{lem17} 
There exist $S_{k,v} \subset M_k^T\mathbb Z^n$
and $a_{k,v,s} \in \mathbb C^*, \, s \in S_{k,v}$ such that
$\zeta_{k,v} := \sum_{s \in S_{k,v}} a_{k,v,s}\, \delta_s$
and $|\zeta_{k,v}|$ are tempered Radon measures and
\begin{equation}
	H_{h,k,v} = \, <\zeta_{k,v}, h>, \ \ h \in \mathcal S_c(\mathbb R^n).
\end{equation}
\end{lem}
Proof. Since $R_{h,k}(\rho_n(\beta_kw)) = 
\sum_{\ell \in \mathbb Z^n} h(\beta_k\ell) e^{-2\pi i \beta_k \ell \cdot w},$ $H_{h, k,v}$ is a linear combination of $h(\beta_k\ell)$
and the coefficient of $h(\beta_k\ell)$ is the coefficient
of $e^{2\pi i \beta_k \ell \cdot w}$ in the trigonometric series expansion of $\frac{\det \left(\partial_{w_i} p_{k,j}\right)}{p_k}$ and $\rho_n(\beta_k w) \in N_{r(k,v)}.$ The frequencies appearing in this series are in the subgroup generated by $\Omega(P).$ This subgroup equals $M^T\mathbb Z^m.$ 
(\ref{LaurentH}) implies
\begin{equation}\label{expansion}
\frac{1}{p_k(w)} = \frac{e^{-2\pi i \beta_k v \cdot w}}{q_{k,v}}
\sum_{m = 0}^\infty  (g_{k,v}(\rho_n(\beta_kw)))^m.
\end{equation}
Since 
$|g_{k,v}(\rho_n(\beta_k w))| \leq \frac{1}{2},$ 
the coefficient of each exponential
$e^{2\pi i \beta_k \ell \cdot w}$ 
in this expansion is bounded by 
$\frac{2}{q_{k,v}}.$ 
The set $A$ of frequencies of 
$g_{k,v}(\rho_n(\beta_kw))$ is
finite and contained in the interior of a closed pointed cone (contains no line) defined by the vertex $v$ of the Newton polytope $\mathcal N(p_k).$ Therefore the set  $B$ of frequencies of the expansion satisfies
\begin{equation}\label{freqexpansion}
	B \subset -v + (\{0\} \cup A \cup (A + A) \cup (A + A + A) \cup \cdots ).
\end{equation}
A combinatorial argument shows that there exists $b > 0$ such that
\begin{equation}\label{boundB}
	|B_n(0,r) \cap B| = O(r^b), \ \ r \rightarrow \infty.
\end{equation}
Multiplication by the trigonometric polynomial
$\det \left(\partial_{w_i} p_{k,j}\right)$ preserves the asymptotic polynomial bound on the distribution of frequencies.
This concludes the proof.
\begin{lem}\label{lem18} 
There exist $S_k \subset M_k^T\mathbb Z^n$
and $a_{k,s} \in \mathbb C^*, \, s \in S_k$ such that
$\zeta_k := \sum_{s \in S_k} a_{k,s}\, \delta_s$
and $|\zeta_k|$ are tempered Radon measures and
\begin{equation}\label{Pkzetak}
	<\widehat {div P_k}, h>\, = \, <\zeta_k, h>, \ \ h \in \mathcal S_c(\mathbb R^n).
\end{equation}
\end{lem}
Proof. Follows from preceding results by defining 
\begin{equation}\label{setak}
\zeta_k := (-1)^n \sum_{v \in \mathcal N(p_k)} C(k,v)\, \zeta_{k,v}
\end{equation}
where $C(k,v)$ are the combinatorial coefficients.
\begin{lem}\label{lem19}  
There exists a discrete subset $S \subset \mathbb R^n$ and 
a sequence $a_s \in \mathbb C^{*}$ such that 
$\zeta := \sum{s \in S} a_s \delta_s$ and $|\zeta|$ are tempered Radon measures and
\begin{equation}\label{convzeta}
\lim_{k \rightarrow \infty} <\zeta_k,h> \, = \, <\zeta, h>, \ \ h \in \mathcal S_c(\mathbb R^n).
\end{equation}
\end{lem}
Proof.
As $k \rightarrow \infty,$ 
$\mathcal N(P_k) \rightarrow \mathcal N(P)$ and the corresponding coefficients converge. Therefore
$B_n(0,r) \cap S_k$ is a Cauchy sequence in the Hausdorff metric for every $r > 0$ so $S_k$ converges pointwise to a discrete set $S$ and the coefficients $a_{k,s}$
converge to $a_s$ so that $\zeta$ defined above satisfy the asserted properties.
\begin{cor}\label{finalcor}
\begin{equation}
	<\widehat {div P}, h>\, = \, <\zeta, h>, \ \ h \in \mathcal S_c(\mathbb R^n).
\end{equation}
\end{cor}
Proof. Proposition \ref{prop8}, (\ref{Pkzetak}) and (\ref{convzeta}) imply 
$$
<\widehat {div P}, h> \, = 
\, \lim_{k \rightarrow \infty} \, <\widehat  {div P_k}, h> \, =
\lim_{k \rightarrow \infty} \, <\zeta_k, h>  \, = \, <\zeta,h>.
$$
This conclude the proof of Proposition \ref{prop10}.
\subsection{Proof of Theorem 3}\label{subsec5.4}
Assume that $P \in \mathcal T_n^n$ is real-rooted and uniformly generic. If $\mathcal N(P)$ is unfolded then Proposition \ref{prop10} implies that $div\, P$ is a Fourier quasicrystal on $\mathbb R^n.$ If $\mathcal N(P)$ is folded then for every $j = 1,...,n$ and 
$v \in \mathcal V(p_j)$ choose a random unit vector $u_v$ and then for every $\epsilon > 0$ define
$p_{j,\epsilon}$ be the trigonometric polynomial $p_j$ whose frequency $v$ replaced by frequency $v+\epsilon u_v$ and whose amplitude at $v+\epsilon u_v$ equals the amplitude of $p_j$ at $v.$
Define 
$P_\epsilon := [p_{1,\epsilon},...,p_{n,\epsilon}]^T.$ 
Since foldedness is defined by equalities, with probability one the trigonometric map $P_\epsilon$ is unfolded.
Therefore there exist unit vectors $u_v$ such that 
$\mathcal N(P_\epsilon)$ is unfolded for every $\epsilon > 0$. Then Proposition \ref{prop10} implies $div\, P_\epsilon$ is a Fourier quasicrystal on $\mathbb R^n.$ Since $\Delta(div\, P) < \infty$
and $\Delta(div\, P_\epsilon)$ are uniformly bounded for sufficiently small $\epsilon,$ the following analogue of Proposition \ref{prop8} holds.
\begin{equation}\label{anprop8}
<\widehat {div P}, h> \, = \, \lim_{\epsilon \rightarrow 0} \, <\widehat  {div P_\epsilon}, h>, \ \ h \in \mathcal S_c(\mathbb R^n).
\end{equation}
$P_\epsilon$ converges to $P$ uniformly on bounded subsets of $\mathbb C^n$ so the same argument used to prove Proposition \ref{prop10} implies that
$div\, P$ is a Fourier quasicrystal.
%
%
\section{Conclusions and Research Questions}\label{sec6}
Section 2 provided a background in harmonic and functional analysis relevant to almost periodic measures and Fourier quasicrystals. Proposition \ref{prop2} implies that every Fourier quasicrystal on $\mathbb R^n$ is a Bohr almost periodic measure. Example 2 constructs a set whose associated measure is a Bohr almost periodic Poisson measure with a dense spectrum. It raises the following question which was asked by Favorov in a personal communication:
\begin{question}\label{question1}
Do there exist multisets in $\mathbb R^n$ whose associated measures are Bohr almost periodic but not Poisson?
\end{question}
If $L \subset \mathbb R^n$ is a lattice subgroup and 
$a : L \mapsto \mathbb R^n$ 
is Bohr almost periodic, then Bochner's result (\ref{Bochner}) implies that the multiset parameterized by the function $s : \mathbb Z \mapsto \mathbb R$ defined by 
$s(\ell) := \ell + a(\ell) : \ell \in L,$ 
is Bohr almost periodic, and hence its associated measure 
$\mu := \sum_{\ell \in L} \delta_{s(\ell)}$ is 
Bohr almost periodic. Favorov (\cite{favorov8}, Theorem 1) proved that all Bohr almost periodic multisubsets in $\mathbb R$ arise in this way. The Fourier series of the continuous 
$$g(x) := 
\sum_{k = 1}^\infty \frac{\sin 2\pi kx}{k}
= \frac{i}{2} \ln \frac{1-e^{-2\pi i x}}{1-e^{2\pi i x}}, 
\ \ x \in \mathbb R/\mathbb Z,$$ 
does not converge absolutely.
Choose $\alpha \in \mathbb R \backslash \mathbb Q$ and define 
$a : \mathbb Z \rightarrow \mathbb R$ by
$$a(\ell) := g(\ell\alpha), \ \ \ell \in \mathbb Z.$$ 
Then $a$ is Bohr almost periodic and its spectrum is a dense subset of $\mathbb R/\mathbb Z.$ Moreover,
its Fourier-Bohr transform $\mathcal F_B(a) : \mathbb R/\mathbb Z \mapsto \mathbb C$ is not even absolutely convergent when restricted to any nonempty open subset of $\mathbb R/\mathbb Z.$
This pathological behavior suggests:
\begin{conj}\label{conjecture1}
The Bohr almost periodic measure $\mu$ associated to the multiset 
parameterized by the function $\ell + a(\ell), \ell \in \mathbb Z,$ is not a Poisson measure.
\end{conj}
$\mu$ is of toral type with compactification $(\mathbb T^2,\psi),$ hence its Fourier-Bohr transform has a expression derived in \cite{lawton}. Therefore the 
criteria for absolute convergence described by Zygmund \cite{zygmund} may be useful to answer this conjecture.
\\ \\
The Fourier quasicrystals constructed in Example 3 and in Section 5 have positive imteger weights and are of toral type. 
This raises the following:
\begin{question}\label{question2}
Are all Fourier quasicrystals with positive imteger weights of toral type?
\end{question}
The result of Olevsly and Ulanovsky \cite{olevskiiulanovskii2} implies that every Fourier quasicrystal on $\mathbb R$ with positive imteger weights is of toral type. Moreover, we have been unable to construct Fourier quasicrystals on $\mathbb R^n, n > 1$ with positive imteger weights that are not of toral type.
\\ \\
Example 3 constructed Fourier quasicrystals using algebraic varieties that are generally not complete intersections, but we used them in Section 5 
because Proposition 
\ref{prop7} required a result in \cite{bushuevatsikh} about complete intersections. This raises:
\begin{question}\label{question3}
Can the sufficient conditions assumed in Section 5 be relaxed. Specifically, can the complete intersection assumption be removed? 
\end{question}
Question \ref{question3} was answered affirmatively  by Alon, Kummer, Kurasov and Vinzant \cite{alon2}. They constructed Delone Fourier quasicrystals:
$$\Lambda := \psi^{-1}(\mathbb T^m \cap V),$$
where $V \subset \mathbb C^{*m}$ is a strict Lee-Yang variety 
(\cite{alon2}, Definition 2.2) and
$\psi = \rho_m \circ M$ where $M \in \mathbb R^{m \times n}$
is a matrix all of whose $n \times n$ minors are positive. If 
$Q : \mathbb C^{*m} \mapsto \mathbb C^k,$ with $k \geq m,$ is a
Laurent map whose zero set equals $V,$ then $Q$ is $M$-stable 
by Definition \ref{stable1} since  the condition on $M$ defines an 
open set of matrices. Therefore Lemma \ref{lemma13} implies that
the construction in \cite{alon2} includes our constructions in Section 5,
but it is more general because it allows $V$ to be a non complete intersection. 
\\ \\
Example 3 and \cite{alon2} constructed irreducible Delone 
sets whose associated measures are Fourier quasicrystals.
\begin{question}\label{question4}
Is there a method to construct irreducible non Delone Fourier quascrystals?
\end{question}
A proof of the following conjecture may be useful to construct non Delone Fourier quasicrystals of toral type.
\begin{conj}\label{conjecture2}
If $(\Lambda, m)$ is an irreducible non Delone Bohr almost periodic multiset in $\mathbb R^n$ of toral type that is associated with a compactification $(\mathbb T^m,\psi)$ of $\mathbb R^n,$ then $\overline {\psi(\Lambda)}$ is the multipicity variety defined by a continous map 
$$f : \mathbb T^{m-n} \mapsto \mathbb T^m$$
that is transverse to the foliation on $\mathbb T^m$ and not injective. The transversality condition on $f$ implies that it induces an injection of fundamental groups
$$
\pi_1(f) : \pi_1(\mathbb T^{m-n}) \mapsto \mathbb \pi_1(T^m).
$$
\end{conj}
The method used in Section 5 in \cite{lawton} to describe the Delone case must be refined to prove Conjecture \ref{conjecture2}.
\\ \\
{\bf Acknowledgments} The authors thank Lior Alon, Sergey Favorov,  Yves Meyer, Alexander Olevskii, and Peter Sarnak for sharing their knowledge about Fourier quasicrystals and Lang's conjecture. 
We thank the reviewer for suggesting revisions that improved the papers accuracy and clarity.
%

\end{document}